\documentstyle[12pt]{article}
\newtheorem{pro}{Proposition}
\newtheorem{theo}{Theorem}
\newtheorem{cor}{Corollary}
\def\ot{\otimes}
\def\ep{\epsilon}
\def\ra{\longrightarrow}
\def\La{\Lambda}
\def\uqn{U_q(\widehat{\hbox{sl}_N})}
\def\uqpn{U_q^\prime(\widehat{\hbox{sl}_N})}
\def\uqnf{U_q(\hbox{sl}_N)}
\newcommand{\bc}[2]{{}_{#1}C_{#2}}

\title{
Trace Construction of a Basis
for the Solution Space of $sl_N$ qKZ Equation}

\author{
Atsushi Nakayashiki\\
\quad\\
LPTHE, Universite Pierre et Marie Curie,\\
Tour 16, 1er etage, 4 place Jussieu,\\
75252, Paris cedex 05, France\\
e-mail: atsushi@lpthe.jussieu.fr\\
and\\
Graduate School of Mathematics,
Kyushu University,\\
Ropponmatsu 4-2-1, Fukuoka 810, Japan\\
e-mail: atsushi@rc.kyushu-u.ac.jp
}

\date{}
\begin{document}
\maketitle
\begin{abstract}
The trace of intertwining operators
over the level one irreducible highest weight
modules of the quantum affine algebra of
type $A_{N-1}$ is studied.
It is proved that the trace function
gives a basis of the solution space of the
qKZ equation at a generic level.
The highest-highest matrix
element of the composition of intertwining 
operators is explicitly calculated.
The integral formula for the trace is
presented.
\end{abstract}

\section{Introduction}
In this paper we study solutions of the quantized 
Knizhnik-Zamolodchikov (qKZ) equation associated with
the quantum group $\uqnf$.
The idea in this paper stems from the
study of solvable lattice models.

The qKZ equation was introduced in \cite{FR}
as the equation satisfied by the
highest-highest matrix element of the intertwining
operators of quantum affine algebra.
For generic values of parameters the set of matrix elements
give a basis of the solution space over the field of
appropriate periodic functions.
The connection matrix of two solutions with different
asymptotics have been calculated from the commutation
relation of intertwining operators.

The solutions of the qKZ equation associated with 
$U_q(\hbox{sl}_2)$ is systematically studied by Tarasov and
Varchenko \cite{TV} (see also references in it).
In \cite{TV} the solutions are described as the multidimensional 
q-hypergeometric integrals.
It is proved that, for generic values of parameters,
the q-hypergeometric solutions give a basis of the solution space
over the field of appropriate periodic functions.
The connection matrix is determined as the
representation of Felder`s elliptic quantum group.

In this paper we propose another description of the basis
of the solution space of the qKZ equation as 
the traces of intertwining operators
of quantum affine algebra.

Let us consider the $\uqnf$ modules
$V_1$,...,$V_n$ and the trigonometric R matrix 
$R_{ij}(z_i/z_j)$ acting on the tensor product 
$V_i \ot V_j$.
The qKZ equation is the q difference equation
for the $V_1\ot\cdots\ot V_n$ valued function 
$f(z_1,\cdots,z_n)$ of the form
\begin{eqnarray}
f(\cdots ,pz_j,\cdots)&=&
R_{jj-1}(pz_j/z_{j-1})\cdots R_{j1}(pz_j/z_1)
(\kappa^{-H})_j
\nonumber
\\
&&
\times
R_{jn}(z_j/z_n)\cdots R_{jj+1}(z_j/z_{j+1})
f(z_1,\cdots,z_n),
\label{hqkz}
\end{eqnarray}
where $\kappa^{-H}=\prod_{i=1}^{N-1}\kappa_i^{-h_i}$,
$h_1,\cdots,h_{N-1}$ is a basis of the Cartan subalgebra
of $sl_N$ and $(\kappa^{-H})_j$ means that $\kappa^{-H}$
acts on $V_j$.
The complex numbers $p$ and $\kappa_i$`s are the parameters
of the equation.
If we write $p=q^{2(k+N)}$ the number $k$ is called level.

Let $\La_i$ $(0\leq i\leq N-1)$ be the fundamental weights of 
$\widehat{sl_N}$. 
We identify $\La_i$ $(1\leq i\leq N-1)$ with 
the fundamental weights of $sl_N$.
In this paper we consider the case where 
all $V_i$ are isomorphic to
the $N$ dimensional irreducible module $V$
with the highest weight $\La_1$ or $\La_{N-1}$.

Let $V(\La_i)$ be the irreducible highest weight
$\uqn$ module with the highest weight $\La_i$ and
$V_\zeta$ the evaluation module of $V$.
Then there exist, up to normalization, unique
intertwining operators $\Phi(\zeta)$ and $\Psi^\ast(\xi)$:
\begin{eqnarray}
&&
\Phi(\zeta):
V(\La_{i+1})
\longrightarrow
V(\La_{i}) \ot V_\zeta,
\quad
\Psi^\ast(\xi):
V_\xi \ot V(\La_{i})
\longrightarrow
V(\La_{i+1}) .
\nonumber
\end{eqnarray}
We extend the index $i$ of $\La_i$ to the set of integers
and read it by modulo $N$.
The operators $\Phi(\zeta)$ and $\Psi^\ast(\xi)$ are
sometimes called of type I and type II respectively \cite{JM}.
The difference between type I and type II 
is in the place where the evaluation module
is. For type I it is on the right of the highest weight module
while for type II it is on the left.

Denote by $D$ the grading operator of the principal gradation of 
$V(\La_i)$ and 
consider the following trace:
\begin{eqnarray}
&&
G(\zeta_1,\cdots,\zeta_m \vert \xi_1,\cdots,\xi_n \vert x,\kappa)
=
\nonumber
\\
&&
F(\zeta \vert \xi \vert x)^{-1}
\sum_{i=0}^{N-1}
\hbox{tr}_{V(\La_i)}
\Big(
x^D\kappa^H
\Phi(\zeta_1)
\cdots
\Phi(\zeta_m)
\Psi^{\ast}(\xi_n)
\cdots
\Psi^{\ast}(\xi_1)
\Big)
\end{eqnarray}
which is a function taking the value in
$\hbox{Hom}_{{\bf C}}(V^{\ot n},V^{\ot m})$.
Here $F(\zeta \vert \xi \vert x)$ is some scalar function
(cf. (\ref{fun2})).
By the commutation relation of the intertwining operators,
the cyclic property of the trace and the functional equation
of $F(\zeta \vert \xi \vert x)$,
this function satisfies the equations:
\begin{eqnarray}
G({\bf \zeta}\vert \cdots x\xi_i\cdots\vert x,y)
&=&
G({\bf \zeta}\vert{\bf \xi}\vert x,\kappa)
\bar{R}_{ii+1}(\xi_i/\xi_{i+1})
\cdots
\bar{R}_{in}(\xi_i/\xi_{n})
(\kappa^{-H})_{\xi_i}
\nonumber
\\
&&
\times
\bar{R}_{i1}(x\xi_i/\xi_1)
\cdots
\bar{R}_{ii-1}(x\xi_i/\xi_{i-1}),
\label{qkz2}
\\
G(\cdots x^{-1}\zeta_i\cdots\vert {\bf \xi} \vert x,\kappa)
&=&
\bar{R}_{jj-1}(x^{-1}\zeta_j/\zeta_{j-1})\cdots 
\bar{R}_{j1}(x^{-1}\zeta_j/\zeta_1)
(\kappa^{-H})_{\zeta_i}
\nonumber
\\
&&
\times
\bar{R}_{jm}(\zeta_j/\zeta_m)\cdots 
\bar{R}_{jj+1}(\zeta_j/\zeta_{j+1})
G({\bf \zeta}\vert{\bf \xi}\vert x,\kappa),
\label{qkz1}
\end{eqnarray}
where $\bar{R}(\zeta)$ is the trigonometric R matrix
(c.f. (\ref{rmatrix})), $\bar{R}_{ii+1}(\xi_i/\xi_{i+1})$
acts non-trivially on $V_{\xi_i}\ot V_{\xi_{i+1}}$ in 
$V^{\ot n}$ etc.
The equation (\ref{qkz1}) has precisely the same form
as the qKZ equation (\ref{hqkz}).
Let
${}^tG$ be the transpose of $G$, that is,
${}^tG \in \hbox{Hom}_{{\bf C}}(V^{\ast\ot m},V^{\ast\ot n})$,
$V^\ast$ being the dual vector space of $V$.
Then, as the equation for ${}^tG$, (\ref{qkz2}) is of the same form
as (\ref{hqkz}).
Since we use the principal gradation in this paper,
to make a precise correspondence between the parameter $x$ and
the parameter $p$ in (\ref{hqkz}) we need to consider $G$ as a function of 
$z_j=\zeta_j^N$ and $u_j=\xi_j^N$.
Then if $x^N=p=q^{2(k+N)}$, ${}^{t}G$ and $G$ satisfy the qKZ equation
of level $k$ and level $-k-2N$ in the variables $u$ and $z$ respectively.

In this paper, if  $x^{-N}=q^{2(k+N)}$, 
we say (\ref{qkz1})
the qKZ equation of level $k$ with the value in $V^{\ot m}$.
%By this definition, if $x^N=q^{2(k+N)}$, $G$ satisfies
%the qKZ equation of level $k$ with the value in $V^{\ast \ot n}$ 
%in the variables $\xi$ and
%the qKZ equation of level $-k-2N$ with the value in $V^{\ot m}$ 
%in the variables $\zeta$.

Let ${\cal S}^n_k$ and ${\cal S}^{\ast n}_k$ be 
the space of meromorphic solutions
of the qKZ equation of level $k$ with the value in
$V^{\ot n}$ and $V^{\ast\ot n}$ respectively 
and ${\cal F}$ the field of $x$ periodic
meromorphic functions in $n$ variables.
Then the function $G$ defines two maps simultaneously:
\begin{eqnarray}
{}^tG(\zeta \vert \cdot \vert x,\kappa)&:&
V^{\ast \ot m} \ot {\cal F}
\longrightarrow {\cal S}^{\ast n}_k,
\label{normal}
\\
G(\cdot \vert \xi \vert x,\kappa)&:&
V^{\ot n} \ot {\cal F}
\longrightarrow {\cal S}^m_{-k-2N}.
\label{dual}
\end{eqnarray}
In (\ref{normal}), $\zeta_1,\cdots,\zeta_m$ are
parameters of the map and
in (\ref{dual}), $\xi_1,\cdots,\xi_n$ are parameters of the map.
We consider the case $n=m$. We assume $\vert x \vert<1$.
We shall prove that, if $x$ and $\kappa$ are generic,
(\ref{normal}) is an isomorphism for the generic
values of $\zeta_1$,...,$\zeta_m$ and
(\ref{dual}) is an isomorphism for the generic values
of $\xi_1$,...,$\xi_n$.
It is proved by showing that the determinant of $G$
does not vanish identically.
We calculate the determinant at $x=0$ where
$G$ reduces to the highest-highest matrix element.
For the level one irreducible module $V(\La_i)$
the matrix elements can be calculated explicitly
without integral.
This is expected because at $q=1$ we have such formula
calculated by using the Frenkel-Kac bosonization of $V(\La_i)$ \cite{FK}.
In the case of $U_q(\widehat{sl_2})$ the formula
is given in \cite{JM}.
For $\uqn$ we carry out the integral of the integral formula
given by the Frenkel-Jing bosonization in \cite{Ko}
in a similar manner to $N=2$ case.

The case $x=q^2$ is relevant to the physical quantities
in solvable lattice models.
In fact at this value of $x$ if we further specialize the
variables $\zeta_i$ and $\xi_j$ appropriately,
the trace functions give correlation functions 
and form factors of the solvable lattice model constructed from the
R-matrix $\bar{R}(\zeta)$.
We have calculated the determinant of $G$ for $N=n=2$ and $x=q^2$
explicitly.
By the $q$ series expansion we checked that det $G$ does not
vanish identically for $n=3$.
We conjecture that the determinant does not vanish identically
at $x=q^2$. This suggests that the trace description can be effective
for the completeness problem of the space of local fields \cite{S}\cite{BBS}.

The bosonization of intertwining operators
makes it possible not only to derive the integral formula for the
matrix elements but also to derive the integral formula
for the trace.
We have given the integral formula.
Therefore the integral formula for the basis 
of the solution space of (\ref{qkz2}) and (\ref{qkz1})
is given.

The plan of this paper is as follows.

In section 2 we give necessary
notations of quantum affine algebra of type $A^{(1)}_{N-1}$.
We introduce the intertwining operators for the level one
integrable modules in the principal picture in section 3.
In section 4 we give the relation between 
principal picture and homogeneous picture.
It serves for translating the results in the references
into principal picture and vice versa.
In section 5 we introduce the trace of intertwining operators
and derive the equations satisfied by them.
The main results and their proof is given in section 6.
In section 7 
we give an example of the concrete expression of the
determinant of the trace of intertwining operators
in the case of $U_q(\widehat{sl_2})$.
In section 8 we give the integral formula for the
matrix element of the intertwining operators.
The integrated formula for the matrix element 
is given in section 9.
In section 10 the integral formula of the
trace of intertwining operators is presented.
In appendix A we refer the integral formula for the
trace of intertwining operators of $U_q(\widehat{sl_2})$
in \cite{JM}, since in this case it is possible 
to simplify the formula a bit. This simplification
is used in the calculation in the example of section 7.
The bosonic expression of the intertwining operators
are reviewed in appendix B.
The list of the expression of the operators in terms of their
normal ordered operators is given in appendix C.
In appendix D a derivation of the integral formula
for the trace of intertwining operators is explained.

\section{Preliminary}
Let $A=(a_{ij})$ be the generalized Cartan matrix of type
$A^{(1)}_{N-1}$, $\{\alpha_i\}_{i=0}^{N-1}$ and 
$\{h_i\}_{i=0}^{N-1}$ 
the set of simple roots and simple coroots 
respectively so that $<\alpha_i,h_j>=a_{ij}$.

The quantum affine algebra $\uqpn$ is the Hopf algebra
generated by $e_i$, $f_i$, $t_i^{\pm 1}$ $(0\leq i\leq N-1)$
with the following defining relations:
$$
t_it_j=t_jt_i,
\quad
t_i^{\pm1}t_i^{\mp1}=1,
\quad
t_i e_j t_i^{-1}=q^{<h_i,\alpha_j>} e_j,
\quad
t_i f_j t_i^{-1}=q^{-<h_i,\alpha_j>} f_j,
$$

$$
[e_i,f_j]=\delta_{ij}
{t_i-t_i^{-1} \over q-q^{-1}},
\quad
\sum_{k=0}^{N-1}
(-1)^ke_i^{(k)}e_j^{(1-a_{ij}-k)}
=
\sum_{k=0}^{1-a_{ij}}
(-1)^k f_i^{(k)}f_j^{(1-a_{ij}-k)}=0 
\quad i\neq j,
$$
where $e^{(k)}=e^k/[k]!$ and similarly for $f^{(k)}$,
$[k]!=[k]\cdots[2][1]$, $[k]=(q^k-q^{-k})/(q-q^{-1})$.

The coproduct $\Delta$ and the antipode $S$ are given by
$\Delta(e_i)=e_i \ot 1+t_i \ot e_i$,
$\Delta(f_i)=f_i \ot t_i^{-1} + 1 \ot f_i$,
$\Delta(t_i)=t_i \ot t_i$ and
$S(e_i)=-t_i^{-1}e_i$,
$S(f_i)=-f_it_i$,
$S(t_i)=t_i^{-1}$.

We extend the algebra $\uqpn$ by adding the element $D$
such that
$$
[D,e_i]=e_i,
\quad
[D,f_i]=-f_i,
\quad
[D,t_i^{\pm1}]=0,
\quad
\Delta(D)=D\ot 1+1\ot D.
$$
The resulting algebra is denoted by $\uqn$.
We say that an element $X\in \uqn$ is of degree $n$ if
$[D,X]=n$.
We denote by $\La_i$ $(0\leq i\leq N-1)$ the fundamental weights
of $\widehat{sl_N}$.
We identify $\La_i$ $(1\leq i\leq N-1)$ with the fundamental weights
of $sl_N$.
We extend the index $i$ of $\La_i$ to any integer
and read it by modulo $N$.

For a highest weight $\uqn$ module $M$ with the highest weight vector
$v$, $D$ defines a grading on $M$ by $D(Xv)=n$ for a degree $n$ element
$X$ in $\uqn$.
The evaluation representation $V_\zeta=\oplus_{j=0}^{N-1}{\bf C}v_j$
of $\uqpn$ associated with the irreducible $\uqnf$ module with the
highest weight $\La_1$ is given by 
\begin{eqnarray}
&&
f_i v_j = \zeta^{-1} \delta_{ij+1} v_{j+1},\quad
e_i v_j = \zeta \delta_{ij} v_{j-1},\quad
t_i v_j = q^{-\delta_{ij} + \delta_{ij+1}}v_j,
\nonumber
\end{eqnarray}
where the index of $v_j$ should be read modulo $N$.
In terms of $\{\La_j\}$ the weight of $v_j$, 
which we denote $\hbox{wt}v_j$, is given by 
$\hbox{wt}v_j=\La_{j+1}-\La_j$.

We denote the binomial coefficient
by $\bc{n}{r}$, that is, $(1+x)^n=\sum_{r=0}^n\bc{n}{r}x^r$.

In this paper two kinds of variables appear, one is
$u$ and $z$, the other is $\xi$ and $\zeta$.
They are always related by the relation $u=\xi^N$ and $z=\zeta^N$
except in Appendix A where $u=-\xi^2$ and $z=\zeta^2$.

\section{Intertwining operators}
In \cite{DO}\cite{Ko} the evaluation
representation, R matrices and intertwining operators
are described in the homogeneous grading.
We shall rewrite them to the principal picture
so that the description is consistent with the $sl_2$ case
in \cite{JM} and that the equations for the trace of intertwining
operators are free from cumbersome factors.

Let $P\bar{R}(\zeta_1/\zeta_2)$ be the 
intertwiner from $V_{\zeta_1}\ot V_{\zeta_2}$
to $V_{\zeta_2}\ot V_{\zeta_1}$ normalized as
$P\bar{R}(\zeta_1/\zeta_2)(v_0 \ot v_0)=v_0 \ot v_0$,
where $P$ is the permutation operator,
$P(v \ot w)=w \ot v$.
We define the components of $\bar{R}(\zeta)$ by
$\bar{R}(\zeta)(v_i\ot v_j)=
\sum_{i^\prime,j^\prime}
\bar{R}(\zeta)^{ij}_{i^\prime j^\prime}
v_{i^\prime}\ot
v_{j^\prime}.$

They are given by (cf. \cite{DO})
\begin{eqnarray}
&&
{\bar R}(\zeta)_{jj}^{jj}=1,\quad
{\bar R}(\zeta)_{jk}^{jk}=b(\zeta)=
{q(1-\zeta^N) \over 1-q^2\zeta^N}\quad (j\neq k),
\nonumber
\\
&&
{\bar R}(\zeta)_{kj}^{jk}=c_{jk}(\zeta)=
{1-q^2 \over 1-q^2\zeta^N}\zeta^{N\theta(k-j)+j-k}\quad
(j\neq k),
\label{rmatrix}
\end{eqnarray}
where $\theta(k)=1$ $(k\geq 0)$, $=0$ (otherwise) and
$0\leq j,k \leq N-1$.

Let $V(\La_i)$ be the irreducible highest weight
$\uqn$ module with the highest weight $\La_i$ and the 
highest weight vector $\vert \La_i>$, 
$V(\La_i)^\ast$ the restricted dual right highest weight
module of $V(\La_i)$ with the highest weight vector
$<\La_i \vert$ such that $<<\La_i \vert, \vert\La_i>>=1$,
where $<,>$ is the dual pairing.
We denote $<<\La_j \vert, X \vert \La_i>>=
<<\La_j \vert X, \vert \La_i>>=
<\La_j \vert X \vert \La_i>$ for any
$X\in \hbox{Hom}_{{\bf C}}(V(\La_i),V(\La_j))$,
where $X$ acts on $V(\La_j)^{\ast}$ from the right.

The type I and type II intertwining operators 
$\Phi^{(i)}(\zeta)$ and $\Psi^{\ast(i)}(\xi)$
are the $\uqpn$ linear operators of the form
\begin{eqnarray}
&&
\Phi^{(i)}(\zeta): V(\La_{i+1})\ra V(\La_{i})\ot V_\zeta,
\Phi^{(i)}(\zeta)=
\sum_{\ep=0}^{N-1}\Phi_{\ep}^{(i)}(\zeta)\ot v_{\ep},
\nonumber
\\
&&
\Psi^{\ast(i)}(\xi): V_\xi\ot V(\La_i)\ra V(\La_{i+1}),
\Psi^{\ast(i)}(\xi)(v_{\mu}\ot u)=\Psi^{\ast(i)}_{\mu}(\xi)u.
\nonumber
\end{eqnarray}
In the second equation $u\in V(\La_i)$.
We normalize them by the condition that
$$
<\La_i \vert \Phi^{(i)}_{i}(\zeta) \vert \La_{i+1}>=1,
\quad
<\La_{i+1} \vert \Psi^{\ast(i)}_{i}(\zeta) \vert \La_{i}>=1.
$$
Under these normalization the operators 
$\Phi^{(i)}(\zeta)$ and
$\Psi^{(i)\ast}(\xi)$ are unique.

We sometimes omit the upper index $(i)$ of
$\Phi^{(i)}(\zeta)$ and $\Psi^{\ast(i)}(\xi)$
for the sake of simplicity.
They satisfy the following commutation relations (\cite{DO}):

\begin{eqnarray}
&&
R(\zeta_1/\zeta_2)
\Phi(\zeta_1)\Phi(\zeta_2)
=
\Phi(\zeta_2)\Phi(\zeta_1),
\label{rpp}
\\
&&
\Psi^\ast(\xi_2)\Psi^\ast(\xi_1)
R^{\ast}(\xi_1/\xi_2)
=
\Psi^\ast(\xi_1)\Psi^\ast(\xi_2),
\label{ppr}
\\
&&
\Phi(\zeta)\Psi^\ast(\xi)
=
\tau(\zeta/\xi)\Psi^\ast(\xi)\Phi(\zeta),
\label{ppt}
\end{eqnarray}
where
\begin{eqnarray}
&&
\tau(\zeta)=
\zeta^{1-N}
{
\theta_{q^{2N}}(q\zeta^N)
\over
\theta_{q^{2N}}(q\zeta^{-N})
}
\nonumber
\end{eqnarray}
and for any complex number $p$ such that $\vert p \vert<1$ we set 
$\theta_p(z)=(z;p)_\infty(pz^{-1};p)_\infty(p;p)_\infty$,
$(z;p)_\infty=\prod_{k=0}^\infty(1-p^kz)$.

The matrices $R(\zeta)$ 
and $R^{\ast}(\zeta)$
is given by 
$R(\zeta)=r(\zeta){\bar R}(\zeta)$ and
$R^{\ast}(\zeta)=r^{\ast}(\zeta){\bar R}(\zeta)$
with
\begin{eqnarray}
&&
r(\zeta)=
\zeta^{-1}
{
(q^{2N}z^{-1};q^{2N})_\infty
(q^{2}z;q^{2N})_\infty
\over
(q^{2N}z;q^{2N})_\infty
(q^{2}z^{-1};q^{2N})_\infty
},
\quad
r^{\ast}(\zeta)=
-\zeta^{-1}
{
(q^{2N}z^{-1};q^{2N})_\infty
(q^{2N-2}z;q^{2N})_\infty
\over
(q^{2N}z;q^{2N})_\infty
(q^{2N-2}z^{-1};q^{2N})_\infty
}.
\nonumber
\end{eqnarray}

In these commutation relations we use the following notation:
for $v_i \ot v_j\in V_{\zeta_1}\ot V_{\zeta_2}$ and
$v_{j^\prime} \ot v_{i^\prime}\in V_{\zeta_2}\ot V_{\zeta_1}$,
the equation $v_i \ot v_j=v_{j^\prime} \ot v_{i^\prime}$
means
$v_i=v_{i^\prime}$ and $v_j=v_{j^\prime}$.
This is for the sake of simplifying the description of the equation.
Thus in terms of components 
(\ref{rpp}), (\ref{ppr}) and (\ref{ppt}) are written as

\begin{eqnarray}
&&
R(\zeta_1/\zeta_2)_{\ep_1\ep_2}^{\ep_1^\prime\ep_2^\prime}
\Phi_{\ep_1^\prime}(\zeta_1)\Phi_{\ep_2^\prime}(\zeta_2)
=
\Phi_{\ep_2}(\zeta_2)\Phi_{\ep_1}(\zeta_1),
\label{crpp}
\\
&&
R^\ast(\xi_1/\xi_2)^{\ep_1\ep_2}_{\ep_1^\prime\ep_2^\prime}
\Psi^\ast_{\ep_2^\prime}(\xi_2)\Psi^\ast_{\ep_1^\prime}(\xi_1)
=
\Psi^\ast_{\ep_1}(\xi_1)\Psi^\ast_{\ep_2}(\xi_2),
\label{cppr}
\\
&&
\Phi_{\ep}(\zeta)\Psi^\ast_{\mu}(\xi)
=
\tau(\zeta/\xi)\Psi^\ast_{\mu}(\xi)\Phi_{\ep}(\zeta).
\label{cppt}
\end{eqnarray}

Let $\sigma$ be the automorphism of $\uqn$ 
induced by the Dynkin diagram automorphism,
$\sigma(e_i)=e_{i+1}$, 
$\sigma(f_i)=f_{i+1}$,
$\sigma(h_i)=h_{i+1}$.
The indices are understood by modulo $N$.
The Dynkin automorphism $\sigma$ induces the linear automorphism
of $V_\zeta$, the linear isomorphism between the left 
highest weight modules $V(\La_i)$ and $V(\La_{i+1})$,
the linear isomorphism between the right highest weight
modules $V(\La_i)^\ast$ and $V(\La_{i+1})^\ast$ by
$\sigma(v_j)=v_{j+1}$, 
$\sigma(\vert \La_i>)=\vert \La_{i+1}>$,
$\sigma(<\La_i \vert) =<\La_{i-1} \vert$ with the properties
$\sigma(Xv)=\sigma(X)\sigma(v)$ and 
$\sigma(v^\ast X)=\sigma(v^\ast)\sigma^{-1}(X)$
for $X\in\uqn$, $v\in V(\La_i)$, $v^\ast \in V(\La_i)^\ast$.

Then the intertwining operators satisfy the following relations
\begin{eqnarray}
&&
\Phi^{(i)}(\zeta)
=
(\sigma \ot \sigma)
\Phi^{(i-1)}(\zeta) \sigma^{-1},
\quad
\Psi^{\ast(i)}(\xi)=
\sigma
\Psi^{\ast(i-1)}(\xi)
(\sigma^{-1}\ot\sigma^{-1}).
\nonumber
\end{eqnarray}
In terms of components these are
\begin{eqnarray}
&&
\Phi^{(i)}_\ep(\zeta)
=
\sigma
\Phi^{(i-1)}_{\ep-1}(\zeta) \sigma^{-1},
\quad
\Psi^{\ast(i)}_\mu(\xi)=
\sigma
\Psi^{\ast(i-1)}_{\mu-1}(\xi)
\sigma^{-1}.
\nonumber
\end{eqnarray}
These relations are proved by checking the intertwining properties
and the normalization conditions of the right hand side
of the equations using the relation
$$
(\sigma \ot\sigma)\Delta
=\Delta\sigma.
$$

The R matrix ${\bar R}(\zeta)$ is also invariant
with respect to $\sigma$:
$$
{\bar R}(\zeta)
^{\sigma(i)\sigma(j)}
_{\sigma(i^\prime)\sigma(j^\prime)}
=
{\bar R}(\zeta)
^{ij}
_{i^\prime j^\prime},
$$
where $\sigma(i)=i+1$ $(0\leq i\leq N-2)$, 
$\sigma(N-1)=0$.

\section{Principal-homogeneous correspondence}
We shall give relations between the intertwining
operators in this paper
and those in \cite{Ko}\cite{DO}.

Let $V^{(h)}_z=\oplus_{j=0}^{N-1}{\bf C}v_j$ be the homogeneous
evaluation representation given by
\begin{eqnarray}
&&
f_i v_j = z^{-\delta_{i0}} \delta_{ij+1} v_{j+1},\quad
e_i v_j = z^{\delta_{i0}} \delta_{ij} v_{j-1},\quad
t_i v_j = q^{-\delta_{ij} + \delta_{ij+1}}v_j,
\nonumber
\end{eqnarray}
The map $V_\zeta \longrightarrow V^{(h)}_z$ given by
$v_i\mapsto v_i\zeta^i$ commutes with the action of
$\uqpn$, where $z=\zeta^N$.

Let $\tilde{\Phi}^{\La_i V}_{\La_{i+1}}(z)$ and
$\tilde{\Phi}^{V^{\ast} \La_{i+1}}_{\La_{i}}(z)$
be the intertwining operators in \cite{Ko}:
$$
\tilde{\Phi}^{\La_i V}_{\La_{i+1}}(z):
V(\La_{i+1})\ra V(\La_{i})\ot V^{(h)}_z,
\quad
\tilde{\Phi}^{V^{\ast} \La_{i+1}}_{\La_{i}}(z)
:V(\La_i)\ra V^{(h)\ast}_z\ot V(\La_{i+1}).
$$
We set
\begin{eqnarray}
&&
\tilde{\Phi}^{h(i)}(z)=
\tilde{\Phi}^{\La_i V}_{\La_{i+1}}(z),
\quad
\tilde{\Phi}^{V^{\ast} \La_{i+1}}_{\La_{i}}(z)=
\sum_{j=0}^{N-1}v_j^{\ast} \ot \tilde{\Psi}^{\ast h(i)}_j(z),
\nonumber
\end{eqnarray}
where $\{v_j^{\ast}\}$ is the dual basis to $\{v_j\}$,
$<v_i,v_j^\ast>=\delta_{ij}$.
Then
\begin{eqnarray}
&&
\Phi^{(i)}_j(\zeta)=
\zeta^{i-j}\tilde{\Phi}^{h(i)}_j(\zeta^N),
\quad
\Psi^{\ast (i)}_j(\zeta)=
\zeta^{j-i}\tilde{\Psi}^{\ast h(i)}_j(\zeta^N),
\nonumber
\end{eqnarray}
where $0\leq i,j\leq N-1$.
We remark that the dual representation $V^\ast$
in $\tilde{\Phi}^{V^{\ast}\La_{i+1}}_{\La_{i}}(z)$ in \cite{Ko}
is with respect to the antipode inverse.

Let $\bar{R}^{(h)}(z)=\bar{R}_{V^{(1)}V^{(1)}}(z)$ be the 
R matrix in \cite{DO}.
Then
$$
\bar{R}(\zeta)^{ij}_{i^\prime j^\prime}
=
\bar{R}^{(h)}(\zeta^N)^{ij}_{i^\prime j^\prime}\zeta^{i-i^\prime}.
$$

\section{Trace of intertwining operators}
In order to appropriately normalize the trace of 
intertwining operators
we first introduce scalar functions which satisfies
some functional equations.
For complex numbers $p_1,\cdots,p_k$ such that 
$\vert p_i \vert <1$ for any $i$, we define
$$
(z;p_1,\cdots,p_k)_\infty=
\prod_{r_1,\cdots,r_k=0}^\infty(1-p^{r_1}\cdots p^{r_k}z).
$$
We set $\{z\}=(z;q^{2N},x^N)_\infty$ and
\begin{equation}
h^{(\sigma)}(z \vert x)=
{
\{q^{1+\sigma}x^Nz^{-1}\}
\{q^{1+\sigma}z\}
\over
\{q^{2N-1+\sigma}x^Nz^{-1}\}
\{q^{2N-1+\sigma}z\}
},
\label{fun1}
\end{equation}
where $\sigma=0,\pm1$.
Let us define
\begin{equation}
\bar{F}({\bf z} \vert {\bf u} \vert x)
=
\prod_{a<b}h^{(+)}({z_b \over z_a} \vert x)
(\prod_{a,b}h^{(0)}({u_a \over z_b} \vert x))^{-1}
\prod_{a<b}h^{(-)}({u_a \over u_b} \vert x),
\label{fun2}
\end{equation}
and
$$
F({\bf \zeta} \vert {\bf \xi} \vert x)
=
{\bar F}({\bf z} \vert {\bf u} \vert x)
\Big(
{
\prod_{a,b}\theta_x(-\xi_b/\zeta_a)
\over
\prod_{a<b}\theta_x(-\zeta_{b}/\zeta_a)
\prod_{a<b}\theta_x(-\xi_a/\xi_{b})
}
\Big)^{N-1}.
$$

The function $F$ satisfies the following equations:
\begin{eqnarray}
&&
F(\cdots,\zeta_{j+1},\zeta_j,\cdots \vert {\bf \xi} \vert x)
=r(\zeta_j/\zeta_{j+1})
F({\bf \zeta} \vert {\bf \xi} \vert x),
\nonumber
\\
&&
F({\bf \zeta} \vert \cdots,\xi_{j+1},\xi_j,\cdots \vert x)
=r^{\ast}(\xi_j/\xi_{j+1})
F({\bf \zeta} \vert {\bf \xi} \vert x),
\nonumber
\\
&&
F(x^{-1}\zeta_1,\cdots,\zeta_m \vert {\bf \xi} \vert x)
=\prod_{j=1}^n \tau(\xi_j/\zeta_1)
F(\zeta_2,\cdots,\zeta_m,\zeta_1 \vert \xi \vert x),
\nonumber
\\
&&
F({\bf \zeta} \vert x\xi_1,\cdots,\xi_n \vert x)
=\prod_{j=1}^m\tau(\xi_1/\zeta_j)
F(\zeta \vert \xi_2,\cdots,\xi_n,\xi_1 \vert x).
\nonumber
\end{eqnarray}

For complex numbers $y_1$,...,$y_{N-1}$ let us set
$y^{\pm H}=\prod_{j=1}^{N-1}y_j^{\pm h_j}$.
Let $x$ be a complex number satisfying $\vert x \vert<1$.
We define the normalized trace function as
\begin{eqnarray}
&&
G^{(i)}({\bf \zeta}\vert {\bf \xi}\vert x,y)=
{
\hbox{tr}_{V(\Lambda_i)}
\big(
x^Dy^{H}
\Phi(\zeta_1)\cdots\Phi(\zeta_m)
\Psi^{\ast}(\xi_n)\cdots\Psi^{\ast}(\xi_1)
\big)
\over 
F({\bf \zeta} \vert {\bf \xi} \vert x)
}.
\label{ntrace1}
\end{eqnarray}
and set
\begin{eqnarray}
&&
G({\bf \zeta} \vert {\bf \xi} \vert x,{\bf y})
=\sum_{i=0}^{N-1}
G^{(i)}({\bf \zeta} \vert {\bf \xi} \vert x,{\bf y}).
\label{sntrace1}
\end{eqnarray}
These functions take the value in
$\hbox{Hom}_{{\bf C}}(V^{\ot n},V^{\ot m})$.
We define the components of $G$ by
$$
G({\bf \zeta}\vert {\bf \xi}\vert x,y)
(v_{\mu_1}\ot\cdots \ot v_{\mu_n})
=
\sum_{\ep_1,\cdots,\ep_m}
G({\bf \zeta}\vert {\bf \xi}\vert x,y)
^{\ep_1,\cdots,\ep_m}
_{\mu_1,\cdots,\mu_n}
v_{\ep_1}\ot\cdots \ot v_{\ep_m}.
$$

By the functional equation of $F$ and the commutation relations
of intertwining operators
the function $G$ satisfies the following system of equations:

\begin{eqnarray}
&&
{\bar R}_{ii+1}(\zeta_i/\zeta_{i+1})
G(\cdots \zeta_i\zeta_{i+1}\cdots\vert {\bf \xi} \vert x,y)
=
G(\cdots \zeta_{i+1}\zeta_i\cdots\vert {\bf \xi} \vert x,y),
\nonumber
\\
&&
G({\bf \zeta}\vert \cdots \xi_i\xi_{i+1}\cdots \vert x,y)
{\bar R}_{ii+1}(\xi_i/\xi_{i+1})
=
G({\bf \zeta}\vert \cdots \xi_{i+1}\xi_{i} \cdots \vert x,y),
\nonumber
\\
&&
G(x^{-1}\zeta_1,\zeta_2,\cdots,\zeta_m\vert {\bf \xi} \vert x,y)
=
(y^{-H})_{\zeta_1}
G(\zeta_2,\cdots,\zeta_m, \zeta_1 \vert {\bf \xi} \vert x,y),
\nonumber
\\
&&
G({\bf \zeta}\vert x\xi_1,\xi_2,\cdots,\xi_n\vert x,y)
=
G({\bf \zeta}\vert \xi_2,\cdots,\xi_n,\xi_1\vert x,y)
(y^{-H})_{\xi_1},
\nonumber
\end{eqnarray}
where ${\bar R}_{ij}(\zeta_i/\zeta_j)$ acts nontrivially
on the component $V_{\zeta_i}\ot V_{\zeta_j}$ in
$V_{\zeta_1}\ot\cdots\ot V_{\zeta_m}$,
${\bar R}_{ij}(\xi_i/\xi_{j})$
acts nontrivially
on the component $V_{\xi_i}\ot V_{\xi_{j}}$ in
$V_{\xi_1}\ot\cdots\ot V_{\xi_n}$ and 
$(y^{-H})_{\zeta_1}$ means that
$y^{-H}$ acts on the component $V_{\zeta_1}$ etc.
In thses equations we use the same notation as in the equation
(\ref{rpp}), (\ref{ppr}) and (\ref{ppt}) to avoid the use of
permutation operators in the equations (see
the comment before (\ref{crpp}), (\ref{cppr}), (\ref{cppt})).

As a consequence of these equations
$G$ satisfies

\begin{eqnarray}
&&
G(\cdots x^{-1}\zeta_i\cdots\vert {\bf \xi}\vert x,y)
=K_i^{(1)}(\zeta_1,\cdots,\zeta_m\vert x,y)
G({\bf \zeta}\vert{\bf \xi}\vert x,y),
\label{level-k-4}
\\
&&
G({\bf \zeta}\vert \cdots x\xi_i\cdots\vert x,y)
=
G({\bf \zeta}\vert{\bf \xi}\vert x,y)
K_i^{(2)}(\xi_1,\cdots,\xi_n\vert x,y)
\label{levelk}
\\
&&
K_i^{(1)}(\zeta_1,\cdots,\zeta_m\vert x,y)
=
\nonumber
\\
&&
{\bar R}_{ii-1}(x^{-1}\zeta_i/\zeta_{i-1})
\cdots
{\bar R}_{i1}(x^{-1}\zeta_i/\zeta_{1})
(y^{-H})_{\zeta_i}
{\bar R}_{im}(\zeta_i/\zeta_{m})
\cdots
{\bar R}_{ii+1}(\zeta_i/\zeta_{i+1}),
\nonumber
\\
&&
K_i^{(2)}(\xi_1,\cdots,\xi_n\vert x,y)
=
\nonumber
\\
&&
{\bar R}_{ii+1}(\xi_i/\xi_{i+1})
\cdots
{\bar R}_{in}(\xi_i/\xi_{n})
(y^{-H})_{\xi_i}
{\bar R}_{i1}(x\xi_i/\xi_{1})
\cdots
{\bar R}_{ii-1}(x\xi_i/\xi_{i-1}).
\nonumber
\end{eqnarray}

Note that
$$
{}^t K_i^{(2)}(\zeta_1,\cdots,\zeta_m\vert x,y)
=
K_i^{(1)}(\zeta_1,\cdots,\zeta_m\vert x^{-1},y).
$$

If we denote $x^N=q^{2(k+N)}$
then the corresponding equations (\ref{level-k-4})
and the transpose of (\ref{levelk}) are the qKZ equations of level $-k-2N$ 
and level $k$ respectively.

{}From the Dynkin symmetry of the intertwining operators,
$G^{(i)}$ and $G$ satisfy the following equations:
\begin{eqnarray}
&&
G^{(i+1)}
({\bf \zeta} \vert {\bf \xi} \vert x,y_1,\cdots,y_{N-1})
^{\sigma(\ep_1),\cdots,\sigma(\ep_m)}
_{\sigma(\mu_1),\cdots,\sigma(\mu_n)}
=
y_1 G^{(i)}
({\bf \zeta} \vert {\bf \xi} \vert x,y_1^{-1}y_2,\cdots,
y_1^{-1}y_{N-1},y_1^{-1})
^{\ep_1,\cdots,\ep_m}
_{\mu_1,\cdots,\mu_n},
\nonumber
\\
&&
G
({\bf \zeta} \vert {\bf \xi} \vert x,y_1,\cdots,y_{N-1})
=
y_1 G
({\bf \zeta} \vert {\bf \xi} \vert x,y_1^{-1}y_2,\cdots,
y_1^{-1}y_{N-1},y_1^{-1}).
\nonumber
\end{eqnarray}

We define 
${\bar G}^{(i)}({\bf \zeta}\vert {\bf \xi}\vert x,y)$
by the similar formula to (\ref{ntrace1}) where 
$F(\zeta \vert \xi \vert x)$ is replaced by
$\bar{F}({\bf z} \vert {\bf u} \vert x)$.
Then define the function
$\bar{G}({\bf \zeta} \vert {\bf \xi} \vert x,{\bf y})$
similarly to (\ref{sntrace1}).

The function ${\bar G}$ satisfies the following system of equations:

\begin{eqnarray}
&&
\big({\zeta_{i+1}\over \zeta_i}\big)^{N-1}
{\bar R}_{ii+1}(\zeta_i/\zeta_{i+1})
\bar{G}(\cdots \zeta_i\zeta_{i+1}\cdots\vert {\bf \xi} \vert x,y)
=
{\bar G}(\cdots \zeta_{i+1}\zeta_i\cdots\vert {\bf \xi} \vert x,y),
\label{sym1}
\\
&&
{\bar G}({\bf \zeta}\vert \cdots \xi_i\xi_{i+1}\cdots \vert x,y)
({\xi_i\over\xi_{i+1}})^{N-1}{\bar R}_{ii+1}(\xi_i/\xi_{i+1})
=
{\bar G}({\bf \zeta}\vert \cdots \xi_{i+1}\xi_{i} \cdots \vert x,y),
\label{sym2}
\\
&&
{\bar G}(x^{-1}\zeta_1,\zeta_2,\cdots,\zeta_m\vert {\bf \xi} \vert x,y)
=
{\bar G}(\zeta_2,\cdots,\zeta_m, \zeta_1 \vert {\bf \xi} \vert x,y)
\prod_{j=1}^n({\zeta_1\over \xi_j})^{N-1},
\nonumber
\\
&&
{\bar G}({\bf \zeta}\vert x\xi_1,\xi_2,\cdots,\xi_n\vert x,y)
=
{\bar G}({\bf \zeta}\vert \xi_2,\cdots,\xi_n,\xi_1\vert x,y)
\prod_{j=1}^n\big({\zeta_j\over \xi_1}\big)^{N-1}.
\nonumber
\end{eqnarray}

For the Dynkin diagram symmetries exactly the same 
equation as $G^{(i)}$ and $G$ holds for $\bar{G}^{(i)}$ and
$\bar{G}$.

Up to now we do not mention to in which sense
the trace (\ref{ntrace1}) exists and to
the validity of the application of the commutation relations
(\ref{rpp}), (\ref{ppr}) and (\ref{ppt}) inside the trace.
By definition the trace (\ref{ntrace1}) exists
as a formal power series in $x$ whose coefficient
is a finite sum of matrix elements of the intertwining 
operator $\Phi(\zeta_1)\cdots\Psi^{\ast}(\xi_1)$.
It is known that the latter matrix element,
which originally defined as a series in $\zeta$ and $\xi$,
are analytically continued to give a meromorphic function
on $({\bf C}^{\ast})^{n+m}$, where ${\bf C}^{\ast}={\bf C}\backslash \{0\}$
is the algebraic torus.
In fact, as we show in section \ref{trf},
the series in $x$ can be summed up explicitly to express
the trace (\ref{ntrace1}) as a meromorphic function
on $({\bf C}^{\ast})^{n+m}$.
Since the commutation relations 
(\ref{rpp}), (\ref{ppr}) and (\ref{ppt})
hold in the sense of analytically continued
matrix element, they are applicable to 
the series expression of the trace in $x$
and hence to the meromorphic expression of the trace.

\section{Completeness of trace function}
In this section we assume $n=m$ and study the determinant
of $\bar{G}$.
Let $(V^{\ot n})_{k_0,\cdots,k_{N-1}}=
\{v\in V^{\ot n} \vert t_iv=q^{k_{i-1}-k_{i}}v\}$ be the 
weight subspace of $V^{\ot n}$
with respect to $\uqnf$.
From the definition of the trace and the intertwining 
operators, $\bar{G}$ commutes with the action of $t_{i}$ 
$(1\leq i \leq N-1)$.
Therefore the determinant of $\bar{G}$ is the product 
of the determinants taken at each weight subspace.
In the following we fix a set of numbers $k_0$,...,$k_{N-1}$.
The determinant always means the determinant taken
at the weight subspace $(V^{\ot n})_{k_0,\cdots,k_{N-1}}$.

\begin{theo}\label{main}
We assume $\vert x \vert <1$.
Then det $\bar{G}(\zeta \vert \xi \vert x,y)$ 
does not vanish identically as a function of $\zeta_i$`s,
$\xi_j$`s, $x,q$ and $y_k$`s.
Moreover det $\bar{G}(\zeta \vert \xi \vert x,1)$ 
does not vanish identically. Here
$y=1$ means that all $y_i=1$.
\end{theo}

We say that $(x,q,y)=(x,q,y_1,\cdots,y_{N-1})$ is generic if 
det $\bar{G}(\zeta \vert \xi \vert x,y)$ does not vanish
identically as a function of $\zeta$`s and $\xi$`s.
For fixed $(x,q,y)$ we say that the set of complex numbers 
$\zeta_1$,...,$\zeta_n$ is generic
if det $\bar{G}(\zeta \vert \xi \vert x,y)$ does not vanish
identically as a function of $\xi$`s.
Similarly 
for fixed $(x,q,y)$ we say that the set of complex numbers 
$\xi_1$,...,$\xi_n$ is generic
if det $\bar{G}(\zeta \vert \xi \vert x,y)$ does not vanish
identically as a function of $\zeta$`s.

\begin{cor}\label{maincor}
Suppose that $(x,q,y)$ is generic.
The map (\ref{normal}) is an isomorphism for
generic values of $\zeta_1$,...,$\zeta_n$
and
the map (\ref{dual}) is an isomorphism for
generic values of $\xi_1$,...,$\xi_n$.
\end{cor}

To prove Theorem \ref{main} we first express the determinant
of $\bar{G}$ using some single function.
Let us set 
$$
f({\bf \zeta} \vert {\bf \xi} \vert x,{\bf y})
=
\bar{G}({\bf \zeta} \vert {\bf \xi} \vert x,{\bf y})
^{0^{k_0}\cdots (N-1)^{k_{N-1}}}
_{0^{k_0}\cdots (N-1)^{k_{N-1}}},
$$
where $(0^{k_0}\cdots (N-1)^{k_{N-1}})$ means the 
$\ep=(\ep_1,\cdots,\ep_{n})$ such that
$\ep_1\leq \cdots\leq \ep_n$ and the number of $i$ in $\ep$
is $k_i$. We call $f$ the extremal component of $\bar{G}$.

{}From (\ref{sym1}) and (\ref{sym2}) we have,
for $l>k$,
\begin{eqnarray}
{\bar G}(\cdots \zeta_i\zeta_{i+1}\cdots \vert {\bf \xi}
\vert x,y)
_{{\bf \mu}}
^{\cdots lk\cdots}
&=&
a^{(3)}_{kl}(\zeta_i/\zeta_{i+1})
{\bar G}(\cdots \zeta_i\zeta_{i+1}\cdots \vert {\bf \xi}
\vert x,y)
_{{\bf \mu}}
^{\cdots kl\cdots}
\nonumber
\\
&&
+
a_1(\zeta_i/\zeta_{i+1})
{\bar G}(\cdots \zeta_{i+1}\zeta_i\cdots \vert {\bf \xi}
\vert x,y)
_{{\bf \mu}}
^{\cdots kl\cdots},
\label{hecke1}
\\
{\bar G}({\bf \zeta} \vert \cdots \xi_i\xi_{i+1}\cdots
\vert x,y)
_{\cdots lk\cdots}
^{{\bf \ep}}
&=&
a^{(4)}_{kl}(\xi_i/\xi_{i+1})
{\bar G}({\bf \zeta} \vert \cdots \xi_i\xi_{i+1}\cdots
\vert x,y)
_{\cdots kl\cdots}
^{{\bf \ep}}
\nonumber
\\
&&
+
a_2(\xi_i/\xi_{i+1})
{\bar G}({\bf \zeta} \vert \cdots \xi_{i+1}\xi_i\cdots
\vert x,y)
_{\cdots kl\cdots}
^{{\bf \ep}},
\label{hecke2}
\end{eqnarray}
where
$$
a_1(\zeta)=
{
q^{-1}\zeta^{N-1}(1-q^2\zeta^N)
\over
1-\zeta^N
},
\quad
a_2(\zeta)=
{
q^{-1}\zeta^{-N+1}(1-q^2\zeta^N)
\over
1-\zeta^N
},
$$
$$
a^{(3)}_{kl}(\zeta)=-
{
q^{-1}\zeta^{N+k-l}(1-q^2)
\over
1-\zeta^N
},
\quad
a^{(4)}_{kl}(\zeta)=-
{
q^{-1}\zeta^{l-k}(1-q^2)
\over
1-\zeta^N
}.
$$
Using these equations it is possible to express
any component of $\bar{G}$ in terms of $f$.
In order to describe this expression precisely
let us introduce the lexicographical order,
on the set of $\ep=(\ep_1,\cdots,\ep_n)$
such that $\sharp\{j \vert \ep_j=i\}=k_i$,
comparing from left to right.
The minimal element is 
$(0^{k_0}\cdots (N-1)^{k_{N-1}})$.
It is convenient to associate 
$M=(M_0,\cdots,M_{N-1})$ with $\ep=(\ep_1,\cdots,\ep_n)$
such that 
$M_i=\{m^{(i)}_1<\cdots< m^{(i)}_{k_i}\}=
\{j \vert \ep_j=i\}$.
We use both notations to specify components of $\bar{G}$.
We denote the minimal element by $M^0=(M^0_0,\cdots,M^0_{N-1})$.
For $M=(M_0,\cdots,M_{N-1})$ we set
$(\zeta_{M_0},\cdots,\zeta_{M_{N-1}})=
(\zeta_{m^{(0)}_1},\cdots,\zeta_{m^{(N-1)}_{k_{N-1}}})$.
Then
\begin{eqnarray}
&&
{\bar G}(
{\bf \zeta}
\vert 
{\bf \xi}
\vert x,y)
^{M}_{L}
=
\sum_{M^{\prime} \leq M}
a^{MM^{\prime}}
{\bf G}
(
\zeta_{M_0^\prime},\cdots,\zeta_{M_{N-1}^\prime} \vert
{\bf \xi}\vert x,y)^{M^0}_L,
\nonumber
\end{eqnarray}
with 
$$
a^{MM}=\prod_{r>l}\prod_{a\in M_r,b\in M_l a<b}
a_1(\zeta_a/\zeta_b).
$$
Similarly
\begin{eqnarray}
&&
{\bar G}(
{\bf \zeta}
\vert 
{\bf \xi}
\vert x,y)
^{M}_{L}
=
\sum_{L^{\prime} \leq L}
b^{LL^{\prime}}
{\bf G}
(
{\bf \zeta} \vert
\xi_{L_0^\prime},\cdots,\xi_{L_{N-1}^\prime} \vert x,y)^{M}_{M^0},
\nonumber
\end{eqnarray}
with 
$$
b^{LL}=\prod_{r>l}\prod_{a\in L_r,b\in L_l a<b}
a_2(\xi_a/\xi_b).
$$

These equations mean that the matrix $\bar{G}$ and the matrix whose
components consist of $f$ 
with permuted variables are connected by the product of two
triangular matrices 
with diagonal components $(a^{MM})_M$ and $(b^{MM})_M$
respectively.
Thus we have

\begin{pro}
$$
\hbox{det}({\bar G}({\bf \zeta}\vert {\bf \xi}\vert x,y)
^{{\bf \epsilon}}
_{{\bf \mu}})=
\Big(
\prod_{i<j}a_1(\zeta_i/\zeta_j)a_2(\xi_i/\xi_j)
\Big)^{n(k_0,\cdots,k_{N-1})}
\qquad\qquad\qquad
$$
$$
\qquad\qquad\qquad\qquad
\times
\hbox{det}
\Big(
f({\bf \zeta}_{M_0},\cdots,{\bf \zeta}_{M_{N-1}} \vert
{\bf \xi}_{L_0},\cdots,{\bf \xi}_{L_{N-1}}\vert x,y)
\Big)_{M,L},
$$
where
$$
n(k_0,\cdots,k_{N-1})=
\sum_{0\leq l< r\leq N-1}
n_{l,r}(k_0,\cdots,k_{N-1}),
$$
$$
n_{l,r}(k_0,\cdots,k_{N-1})=
\prod_{j=0}^{N-2}
\bc{n-2-\sum_{i=0}^{j-1}k_i^\prime}{k_j^\prime},
$$
$$
(k_0^\prime,\cdots,k_{N-1}^\prime)
=
(k_0,\cdots,k_l-1, k_{l+1},\cdots,k_{r-1}, k_r-1,\cdots,k_{N-1}),
$$
where the empty sum $\sum_{i=0}^{-1}k_i^\prime$ should be
understood as $0$.
\end{pro}

Note that, by definition, 
${\bar G}({\bf \zeta}\vert {\bf \xi}\vert x,y)$
reduces to the matrix element at $x=0$:
\begin{eqnarray}
{\bar G}({\bf \zeta}\vert {\bf \xi}\vert 0,y)
^{{\bf \epsilon}}
_{{\bf \mu}}&=&
\bar{F}(z \vert u\vert 0)
\nonumber
\\
&&
\times
\sum_{i=0}^{N-1}
<\Lambda_i\vert
y^H
\Phi_{\ep_1}(\zeta_1) \cdots \Phi_{\ep_m}(\zeta_m)
\Psi^{\ast}_{\mu_n}(\xi_n) \cdots \Psi^{\ast}_{\mu_1}(\xi_{1})
\vert\Lambda_i>.
\nonumber
\end{eqnarray}
The matrix element can be calculated explicitly using the
bosonization of the intertwining operators \cite{Ko}
on the Frenkel-Jing bosonization of $V(\La_i)$ \cite{FJ}.
In fact the bosonization gives the integral formula
for the matrix elements.
The integral of the extremal component can be calculated 
rather easily.
For $N=2$ such calculation is done in \cite{JM}.
We give the integral formula in section \ref{matint}
and the integrated formula with its derivation
in section \ref{matrat}.

By specializing the formula of Proposition \ref{integrated}
in section \ref{matrat} to $i=j$, $k_r=l_r$ 
$(0\leq r\leq N-1)$, multiplying $y_i$ and 
summing up in $i$ from $0$ to $N-1$
we get

\begin{pro}\label{extrem}
We have
$$
f({\bf \zeta} \vert {\bf \xi} \vert 0,y)
=
C
\prod_{a=1}^n
({\zeta_a \over \xi_a})
^{
(N-1)(1-a)
}
\prod_{r=0}^{N-1}
\prod_{a\in M^0_r}
({\xi_a \over\zeta_a})^r
\quad\quad\quad\quad\quad\quad
$$
$$
\quad\quad\quad\quad\quad
\quad\quad\quad\quad\quad
\times
\prod_{r<l}
\prod_{a\in M^0_r, b\in M^0_l}
{
(z_a-qu_b)(u_a-qz_b)
\over
(z_a-q^2z_b)(u_a-q^2u_b)
}
\sum_{i=0}^{N-1}
y_i
\prod_{a=1}^n
({\zeta_a \over \xi_a})^i
\prod_{a=1}^{K_{i-1}}
({u_a \over z_a})
$$
where $K_j=k_0+\cdots+k_{j}$, $K_{-1}=0$, $y_0=1$,
$$
C
=
(-1)^{
\sum_{r=0}^{N-2}(r+1)k_r
}
q^{
{1\over2}K_{N-1}^2
-{1\over2}\sum_{r=0}^{N-2}k_r^2
+K_{N-2}k_{N-1}
}.
$$
The empty product from $1$ to $0$ should be understood as one.
\end{pro}
\vskip5mm

\noindent
Proof of Theorem \ref{main}.
It is sufficient to prove that
$\det \Big(
f({\bf \zeta}_{M_0},\cdots,{\bf \zeta}_{M_{N-1}} \vert
{\bf \xi}_{L_0},\cdots,{\bf \xi}_{L_{N-1}}\vert 0,1)
\Big)_{M,L}$ does not vanish identically, where
$1$ means $y_i=1$ for any $i$.
Let $P=\prod_{a=0}^{N-2}\bc{n-K_{a-1}}{k_a}$ be 
the size of the determinant.
We set
$$
f_1(\zeta_1,\cdots,\zeta_n)=
\prod_{a=1}^n
\zeta_a^{(N-1)(1-a)}
\prod_{r=0}^{N-1}\prod_{a\in M_r^0}\zeta_a^{-r},
\quad
f_2(\xi_1,\cdots,\xi_n)=f_1(\xi_1,\cdots,\xi_n)^{-1}.
$$
Then by Proposition \ref{extrem} we have
\begin{eqnarray}
&&
\det \Big(
f({\bf \zeta}_{M_0},\cdots,{\bf \zeta}_{M_{N-1}} \vert
{\bf \xi}_{L_0},\cdots,{\bf \xi}_{L_{N-1}}\vert 0,1)
\Big)_{M,L}
\nonumber
\\
&&
=
C^P\prod_M
f_1(\zeta_{M_0},\cdots,\zeta_{M_{N-1}})
f_2(\xi_{M_0},\cdots,\xi_{M_{N-1}})
\nonumber
\\
&&
\times
\prod_M\prod_{r<l}\prod_{a\in M_r,b\in M_l}
(z_a-q^2z_b)^{-1}(u_a-q^2u_b)^{-1}
{\cal D}(\zeta \vert\xi \vert q),
\nonumber
\end{eqnarray}
where ${\cal D}(\zeta \vert\xi \vert q)
=\hbox{det} \big({\cal D}(\zeta \vert\xi \vert q)_{M,L}\big)_{M,L}$ and
\begin{eqnarray}
&&
{\cal D}(\zeta \vert\xi \vert q)_{M,L}
\nonumber
\\
&&
=
\prod_{r<l}
\Big[
\prod_{a\in M_r, b\in L_l}
(z_a-qu_b)
\prod_{a\in L_r, b\in M_l}
(u_a-qz_b)
\Big]
\sum_{j=0}^{N-1}
\prod_{a=1}^n
({\zeta_a \over \xi_a})^j
\prod_{r=0}^{j-1}
\big(
\prod_{a\in M_r}z_a^{-1}
\prod_{a\in L_r}u_a
\big).
\nonumber
\end{eqnarray}
We consider the case $q=1$ and $\zeta_j=\xi_j$ for any $j$.
It is easy to see that if $M$ is different from $L$ then 
${\cal D}(\zeta \vert\xi \vert 1)_{M,L}=0$.
Hence the matrix 
$\big({\cal D}(\zeta \vert\zeta \vert 1)_{M,L}\big)_{M,L}$
is a diagonal matrix and 
$$
\hbox{det}\big({\cal D}(\zeta \vert\zeta \vert 1)_{M,L}\big)_{M,L}
=
\prod_{M}{\cal D}(\zeta \vert\zeta \vert 1)_{M,M}
=
N^P
\prod_M
\prod_{r<l}
\prod_{a\in M_r,b\in M_l}
(z_a-z_b)^2.
$$
This completes the proof.

\section{Determinant formula at $N=2$ and $x=q^2$}
Let us give here examples of the
explicit formulae for the determinants of ${\bar G}$.
To understand the general structure of the formula in the example below
we first present the system of equations satisfied by the determinant
and give one solution of it.
By taking the determinant of the equation of $\bar{G}$
at the weight space $(V^{\ot n})_{k_0,k_1}$
we have
\begin{eqnarray}
&&
({\zeta_{i+1}\over \zeta_i})^{\bc{n}{k_0}}
\Big({
z_i-q^2z_{i+1}
\over
z_{i+1}-q^2z_i
}\Big)^{\bc{n-2}{k_0-1}}
\hbox{det}{\bar G}(\cdots \zeta_i\zeta_{i+1}\cdots\vert {\bf \xi} \vert x,y)
=
\hbox{det}{\bar G}(\cdots \zeta_{i+1}\zeta_i\cdots\vert {\bf \xi} \vert x,y),
\nonumber
\\
&&
({\xi_i\over\xi_{i+1}})^{\bc{n}{k_0}}
\Big({
u_i-q^2u_{i+1}
\over
u_{i+1}-q^2u_i
}\Big)^{\bc{n-2}{k_0-1}}
\hbox{det}{\bar G}({\bf \zeta}\vert \cdots \xi_i\xi_{i+1}\cdots \vert x,y)
=
\hbox{det}{\bar G}({\bf \zeta}\vert \cdots \xi_{i+1}\xi_{i} \cdots \vert x,y),
\nonumber
\\
&&
\hbox{det}{\bar G}
(x^{-1}\zeta_1,\zeta_2,\cdots,\zeta_n\vert {\bf \xi} \vert x,y)
=
\hbox{det}{\bar G}
(\zeta_2,\cdots,\zeta_n, \zeta_1 \vert {\bf \xi} \vert x,y)
(-1)^{(n-1)\bc{n-2}{k_0-1}}
(\prod_{j=1}^n{\zeta_1\over \xi_j})^{\bc{n}{k_0}},
\nonumber
\\
&&
\hbox{det}{\bar G}
({\bf \zeta}\vert x\xi_1,\xi_2,\cdots,\xi_n\vert x,y)
=
\hbox{det}{\bar G}
({\bf \zeta}\vert \xi_2,\cdots,\xi_n,\xi_1\vert x,y)
(-1)^{(n-1)\bc{n-2}{k_0-1}}
(\prod_{j=1}^n{\zeta_j\over \xi_1})^{\bc{n}{k_0}},
\nonumber
\end{eqnarray}

If $x=q^2$ and $y=1$, one solution to these system of equations is
given by
\begin{eqnarray}
&&
Q({\bf \zeta}\vert {\bf \xi})=
\nonumber
\\
&&
\Big(
\prod_{j=1}^n({\xi_j\over\zeta_j})^{n-1}
\prod_{k<k^\prime}
{
z_{k^\prime}-q^2z_k
\over
u_k-q^2u_{k^\prime}
}
\Big)^{\bc{n-2}{k_0-1}}
\Big(
\prod_{j=1}^n({\xi_j\over\zeta_j})^{j-1}
\theta_{q^2}(-
\prod_{j=1}^n{\xi_j\over\zeta_j})
\Big)^{\bc{n}{k_0}}.
\nonumber
\end{eqnarray}
Any other meromorphic solution of the equation
is given by multiplying a meromorphic function which is
symmetric and $q^2$ periodic in $\zeta_i$`s and $\xi_i$`s
respectively to $Q({\bf \zeta}\vert {\bf \xi})$.

\vskip5mm

\noindent
{\bf Example 1.}
We consider the case of $k_0=0$.
The formula is from \cite{JM}.
\begin{eqnarray}
&&
{\bar G}(\zeta_1,\cdots,\zeta_n \vert \xi_1,\cdots,\xi_n
\vert x,y)
_{-\cdots-}^{-\cdots-}
=
(x^2)_\infty
\prod_{j=1}^n\Big({\xi_j\over\zeta_j}\Big)^j
\theta_x(-y\prod_{j=1}^n{\zeta_j\over\xi_j}).
\nonumber
\end{eqnarray}
\vskip5mm

\noindent
{\bf Example 2.}
Let us consider the case $x=q^2$ $y=1$, $n=2$, $k_0=1$.
Then
\begin{eqnarray}
&&
{\bar G}(\zeta_1,\zeta_2 \vert \xi_1,\xi_2\vert q^2,1)^{+-}_{+-}
=
\nonumber
\\
&&
q(q^4)_\infty
\prod_{j=1}^2({\xi_j\over\zeta_j})
{
\theta_{q^2}(-q\prod_{j=1}^2{\xi_j\over\zeta_j})
\over 
u_1-q^2u_2
}u_2
(1-{\zeta_1\xi_1 \over\zeta_2\xi_2}),
\end{eqnarray}
and
\begin{eqnarray}
&&
\hbox{det}\big({\bar G}(\zeta_1,\zeta_2 \vert \xi_1,\xi_2)
^{\epsilon_1\epsilon_2}_{\mu_1\mu_2}\big)
=
\nonumber
\\
&&
-q^{2}(q^4)_{\infty}^2
\prod_{j=1}^2({\xi_j\over\zeta_j})^{2j}
{z_2-q^2z_1 \over u_1-q^2u_2}
\theta_{q^2}(-q\prod_{j=1}^2{\xi_j\over\zeta_j})^2.
\end{eqnarray}
Here $(z)_\infty=(z:q^4)_\infty$ and 
$z_i=\zeta_i^2$, $u_i=\xi_i^2$.
This formula is calculated using the technique
found in \cite{NPT}.

By using (\ref{hecke1}), (\ref{hecke2}) we can calculate
other components of $\bar{G}$.
By Dynkin diagram symmetry we know a priori
that
$\bar{G}(\zeta\vert\xi\vert q^2,1)^{-+}_{-+}=
\bar{G}(\zeta\vert\xi\vert q^2,1)^{+-}_{+-}$,
$\bar{G}(\zeta\vert\xi\vert q^2,1)^{-+}_{+-}=
\bar{G}(\zeta\vert\xi\vert q^2,1)^{+-}_{-+}$.
Then the concrete expression for $\bar{G}$ is given by
\begin{eqnarray}
&&
\bar{G}(\zeta\vert\xi\vert q^2,1)(v_{+}\ot v_{-})
=
\nonumber
\\
&&
(q^4)_\infty
\prod_{j=1}^2({\xi_j\over\zeta_j})
{
\theta_{q^2}(-q\prod_{j=1}^2{\xi_j\over\zeta_j})
\over 
u_1-q^2u_2
}u_2
\Big(
q(1-{\zeta_1\xi_1 \over\zeta_2\xi_2})v_{+}\ot v_{-}
-{\xi_1 \over\xi_2}(1-q^2{\zeta_1\xi_2 \over \zeta_2\xi_1})
v_{-}\ot v_{+}
\Big),
\nonumber
\\
&&
\bar{G}(\zeta\vert\xi\vert q^2,1)(v_{-}\ot v_{+})
=
\nonumber
\\
&&
(q^4)_\infty
\prod_{j=1}^2({\xi_j\over\zeta_j})
{
\theta_{q^2}(-q\prod_{j=1}^2{\xi_j\over\zeta_j})
\over 
u_1-q^2u_2
}u_2
\Big(
-{\xi_1 \over\xi_2}(1-q^2{\zeta_1\xi_2 \over \zeta_2\xi_1})
v_{+}\ot v_{-}
+
q(1-{\zeta_1\xi_1 \over\zeta_2\xi_2})
v_{-}\ot v_{+}
\Big).
\nonumber
\end{eqnarray}

\section{Integral formula for matrix elements}\label{matint}
We set
$$
\bar{G}^{(ij)}({\bf \zeta} \vert {\bf \xi})
^{\ep_1,\cdots,\ep_m}
_{\mu_1,\cdots,\mu_n}
=
{
<\La_i \vert
\Phi_{\ep_1}(\zeta_1)
\cdots
\Phi_{\ep_{m}}(\zeta_m)
\Psi^{\ast}_{\mu_n}(\xi_n)
\cdots
\Psi^{\ast}_{\mu_1}(\xi_1)
\vert \La_j>
\over
\bar{F}(z \vert u\vert 0)
},
$$
where $\bar{F}(z \vert u\vert 0)$ is given by (\ref{fun1}) 
and (\ref{fun2}).
We need to assume $j+n=i+m$ mod.$N$ for 
the matrix element to be well defined.

Let $k_r=\sharp\{j \vert \ep_j=r\}$,
$l_r=\sharp\{j \vert \mu_j=r\}$ for $0\leq j\leq N-1$.
Then $m=\sum_{r=0}^{N-1}k_r$ and
$n=\sum_{r=0}^{N-1}l_r$.
The function $\bar{G}^{(ij)}({\bf \zeta} \vert {\bf \xi})
^{{\bf \ep}}_{{\bf \mu}}$ is zero unless
$$
\sum_{r=1}^n\hbox{wt}v_{\mu_r}+\La_j
=
\sum_{r=1}^m\hbox{wt}v_{\ep_r}+\La_i.
$$
Since $\hbox{wt}v_r=\La_{r+1}-\La_r$ 
this condition is written as
\begin{eqnarray}
&&
k_r-l_r=k_{r-1}-l_{r-1}+\delta_{r,i}-\delta_{r,j}
\qquad
0\leq r \leq N-1,
\label{wcond}
\end{eqnarray}
where we understand $k_{-1}=k_{N-1}$ and $l_{-1}=l_{N-1}$.
In particular we have $m-n=j-i+Nr_0$,
where $r_0=k_{N-1}-l_{N-1}$.
We assume the condition (\ref{wcond}).
We set $w^{(a)}_N=q^{N+1}z_a$ and $v^{(b)}_N=q^{N+1}u_b$
for the sake of convenience.
Then

$$
\bar{G}^{(ij)}({\bf \zeta} \vert {\bf \xi})
^{\ep_1,\cdots,\ep_m}
_{\mu_1,\cdots,\mu_n}
\qquad\qquad\qquad\qquad\qquad\qquad
\qquad\qquad\qquad\qquad\qquad\qquad\qquad
\qquad\qquad\qquad\qquad
$$

$$
=
\bar{C}^{(ij)}(\ep,\mu)
\prod_{a=1}^m
\zeta_a^{(N-1)(m-n+1-a)+j-\ep_a
}
\prod_{b=1}^n
\xi_b^{
(N-1)(b-1)-j+\mu_b
}
\qquad\qquad\qquad\qquad\qquad\qquad\qquad\qquad
$$

$$
\times
\int_{C^{(1)}_{\ep_1+1}}{dw^{(1)}_{\ep_1+1} \over 2\pi i}
\cdots
\int_{\tilde{C}^{(n)}_{N-1}}{dv^{(n)}_{N-1} \over 2\pi i}
\prod_{a;\ep_a\leq j-1}(w^{(a)}_j)^{-1}
\prod_{b;\mu_b\leq j-1} v^{(b)}_j
\qquad\qquad\qquad\qquad\qquad\qquad\qquad\quad
$$

$$
\times
\prod_{a=1}^m
\prod_{k=\ep_a+1}^{N-1}
{
(q^{-1}-q)w^{(a)}_{k}
\over
(w^{(a)}_k-q^{-1}w^{(a)}_{k+1})
(w^{(a)}_k-qw^{(a)}_{k+1})
}
\prod_{b=1}^n
\prod_{k=\mu_b+1}^{N-1}
{
(q^{-1}-q)v^{(b)}_{k+1}
\over
(v^{(b)}_{k}-q^{-1}v^{(b)}_{k+1})
(v^{(b)}_k-qv^{(b)}_{k+1})
}
\qquad\qquad
$$

$$
\times
\prod_{a<b}
\prod_{k}
{-1\over w^{(a)}_k-qw^{(b)}_{k-1}}
\prod_{a<b}
\prod_{k}
{1\over w^{(a)}_k-qw^{(b)}_{k+1}}
\prod_{a<b}
\prod_{k\leq N-1}
(w^{(a)}_k-q^2w^{(b)}_k)
(w^{(a)}_k-w^{(b)}_k)
\qquad\qquad\qquad
$$

$$
\times
\prod_{a<b}
\prod_{k}
{-1 \over v^{(b)}_{k}-q^{-1}v^{(a)}_{k-1}}
\prod_{a<b}
\prod_{k}
{1 \over v^{(b)}_{k}-q^{-1}v^{(a)}_{k+1}}
\prod_{a<b}
\prod_{k\leq N-1}
(v^{(b)}_k-q^{-2}v^{(a)}_k)
(v^{(b)}_k-v^{(a)}_k)
\qquad\qquad\quad
$$

$$
\times
\prod_{a,b,k}
(w^{(a)}_k-v^{(b)}_{k+1})
\prod_{a,b,k}
(v^{(b)}_{k-1}-w^{(a)}_k)
\prod_{a,b}\prod_{k\leq N-1}
{1 
\over
(w^{(a)}_k-qv^{(b)}_k)
(w^{(a)}_k-q^{-1}v^{(b)}_k)
},
\qquad\qquad\qquad
$$
where
\begin{eqnarray}
&&
\bar{C}^{(ij)}(\ep,\mu)
=
(-1)^{
{1\over2}\sum_{a=1}^m
(N-N\{{i-1+a \over N}\})(N-1-N\{{i-1+a \over N}\})
+
{1\over2}\sum_{b=1}^n
(N-N\{{j-1+b \over N}\})(N-1-N\{{j-1+b \over N}\})
}
\nonumber
\\
&&
\times
(-1)^{-ir_0(N-1)+\delta_{j0}(n+m)(N-1)
+({1\over2}N(N+1)+1)(k_0+l_0)
}
\nonumber
\\
&&
\times
(-1)^{
{1\over2}(j-i)(N-j)(N-j-1)\theta(1\leq i<j)
+{1\over2}(i-j)(N-i)(N-1+i-2j)\theta(i>j\geq 1)
}
\nonumber
\\
&&
\times
(-1)^{
{1\over2}(k_0+l_0)\sum_{r=1}^{N-1}(N-1-r)(N+2+r)(k_r+l_r)
}
\nonumber
\\
&&
\times
q^{{1\over2}(N+1)((i-j)(i-j-1)+r_0^2N(N-1)+2jr_0N-2ir_0(N-1))}.
\nonumber
\end{eqnarray}
Here, for a rational number $r$,
we denote by $\{r\}$ the fractional part of $r$, that is,
$\{r\}=r-[r]$, $[r]$ being the Gauss symbol.
This notation appears only in the description of 
the sign and should not be confused with the
double infinite product.

The integral variables are
$w^{(a)}_k$ $a=1,\cdots,m$, $k=\ep_a+1,\cdots,N-1$
and
$v^{(b)}_k$ $b=1,\cdots,n$, $k=\ep_b+1,\cdots,N-1$.

Each product in $a,b,k$ which appears in the integrand
is over all possible values satisfying the
conditions written in the product symbol.
We must be careful if $w^{(a)}_N$ or $v^{(b)}_N$
appears in the product.
For example
\begin{eqnarray}
&&
\prod_{a<b}
\prod_{k}
{1\over w^{(a)}_k-qw^{(b)}_{k+1}}
\nonumber
\\
&&
=\prod_{1\leq a<b\leq m-k_{N-1}}
\prod_{k=\hbox{max}(\ep_a+1,\ep_b)}^{N-2}
{1\over w^{(a)}_k-qw^{(b)}_{k+1}}
\prod_{1\leq a<b\leq m-k_{N-1}}
{1\over w^{(a)}_{N-1}-qw^{(b)}_{N}}
\nonumber
\\
&&
\times
\prod_{a=1}^{m-k_{N-1}}
\prod_{b=m-k_{N-1}+1}^{m}
{1\over w^{(a)}_{N-1}-qw^{(b)}_{N}}.
\nonumber
\end{eqnarray}
This is because the index $a$ of $w^{(a)}_k$ runs until
$m-k_{N-1}$ if $k\leq N-1$, while $a$ can run until $m$ if
$k=N$.

The integration contours $C^{(a)}_k$ of $w^{(a)}_k$
and $\tilde{C}^{(b)}_k$ of $v^{(b)}_k$ are as follows.

The contour $C^{(a)}_k$ is a simple closed curve
going round the origin in the anticlockwise direction such that
$qw^{(a)}_{k\pm1}$, $q^{\pm1}w^{(b)}_{k\pm1}$ $(a<b)$,
$q^{\pm1}v^{(b)}_{k}$ (any $b$) are inside,
$q^{-1}w^{(a)}_{k\pm1}$ and $q^{\pm1}w^{(b)}_{k\pm1}$ $(a>b)$ are outside.

The contour $\tilde{C}^{(b)}_k$ is a simple closed curve
going round the origin in the anticlockwise direction such that
$q^{-1}v^{(b)}_{k\pm1}$, $q^{\pm1}v^{(a)}_{k\pm1}$ $(a<b)$ are inside,
$qv^{(b)}_{k\pm1}$,
$q^{\pm1}w^{(a)}_{k}$ (any $a$)
and $q^{\pm1}v^{(a)}_{k\pm1}$ $(a>b)$ are outside.

\section{Integrated formula for the extremal component}
\label{matrat}
In this section we give the formula without integration
for the extremal component of the matrix element.
For $0\leq r\leq N-1$ we define $K_r$, $L_r$ by
$$
K_r=\sum_{r^\prime=0}^{r}k_{r^\prime},
\quad
L_r=\sum_{r^\prime=0}^{r}l_{r^\prime}.
$$
We set $K_r=L_r=0$ for $r<0$ or $r\geq N$.
For a proposition $P$ we define
$\theta(P)=1$ if $P$ is true and 
$\theta(P)=0$ otherwise.
The variables are related by
$z_r=\zeta_r^N$, 
$u_r=\xi_r^N$.

\begin{pro}\label{integrated}
We have
$$
\bar{G}^{(ij)}({\bf \zeta} \vert {\bf \xi})
^{0^{k_0}\cdots (N-1)^{k_{N-1}}}
_{0^{l_0}\cdots (N-1)^{l_{N-1}}}
\qquad\qquad\qquad\qquad\qquad\qquad\qquad\qquad\qquad\qquad
$$
$$
=
C^{(ij)}({\bf k} \vert {\bf l})
\prod_{a=1}^m
\zeta_a^{
(N-1)(m-n+1-a)-\ep_a+j
}
\prod_{b=1}^n
\xi_b^{
(N-1)(b-1)+\mu_b-j
}
\prod_{a=1}^{K_{j-1}}z_a^{-1}
\prod_{b=1}^{L_{j-1}}u_b
$$
$$
\times
{
\prod_{a=1}^m
\prod_{b=1}^n
(z_a-qu_b)^{\theta(\ep_a<\mu_b)}
(u_b-qz_a)^{\theta(\ep_a>\mu_b)}
\over
\prod_{a,b=1}^m
(z_a-q^2z_b)^{\theta(\ep_a<\ep_b)}
\prod_{a,b=1}^n
(u_a-q^2u_b)^{\theta(\mu_a<\mu_b)}
},
\qquad\qquad\qquad
$$
where $C^{(ij)}({\bf k} \vert {\bf l})$ is a constant given by
\begin{eqnarray}
&&
C^{(ij)}({\bf k} \vert {\bf l})
\nonumber
\\
&&
=
(-1)^{
{1\over2}\sum_{a=1}^m
(N-N\{{i-1+a \over N}\})(N-1-N\{{i-1+a \over N}\})
+
{1\over2}\sum_{b=1}^n
(N-N\{{j-1+b \over N}\})(N-1-N\{{j-1+b \over N}\})
}
\nonumber
\\
&&
\times
(-1)^{
-ir_0(N-1)+\sum_{r=0}^{N-2}(N-1-r)k_r
+\delta_{j0}(n+m)(N-1)
+({1\over2}N(N+1)+1)(k_0+l_0)
}
\nonumber
\\
&&
\times
(-1)^{
{1\over2}(j-i)(N-j)(N-j-1)\theta(1\leq i<j)
+{1\over2}(i-j)(N-i)(N-1+i-2j)\theta(i>j\geq 1)
}
\nonumber
\\
&&
\times
(-1)^{
{1\over2}(k_0+l_0)\sum_{r=1}^{N-1}(N-1-r)(N+2+r)(k_r+l_r)
}
\nonumber
\\
&&
\times
(-1)^{
\sum_{1\leq a<b \leq K_{N-2}}(N-1-\ep_b)
+
\sum_{1\leq a<b \leq L_{N-2}}(N-1-\mu_b)
+
\sum_{a=1}^{K_{N-2}}
\sum_{b=1}^{L_{N-2}}
(N-1-\hbox{max}(\ep_a,\mu_b))
}
\nonumber
\\
&&
\times
q^{
{1\over2}(N+1)
(
(i-j)(i-j-1)+r_0^2N(N-1)+2jr_0N-2ir_0(N-1)
)
+(j+1)(-K_{j-1}+L_{j-1})
-N\bc{K_{N-2}}{2}
-(N-1)\bc{L_{N-2}}{2}
}
\nonumber
\\
&&
\times
q^{
N(-K_{N-2}+L_{N-2})(k_{N-1}-l_{N-1})
+L_{N-2}l_{N-1}
-\sum_{r=0}^{N-2}(N-1-r)k_r
+\sum_{r=0}^{N-2}(r+1)\bc{k_r}{2}
+\sum_{r=0}^{N-2}r\bc{l_r}{2}
}
\nonumber
\\
&&
\times
q^{
-\sum_{r=0}^{N-2}(r+1)k_rl_r
+NK_{N-2}L_{N-2}
}.
\end{eqnarray}
and $\ep=(0^{k_0},\cdots, (N-1)^{k_{N-1}})$,
$\mu=(0^{l_0},\cdots, (N-1)^{l_{N-1}})$.

\end{pro}

Let us explain how to derive this formula.
First we carry out the integration in the variable $w$ by the order
$$
w^{(1)}_1\rightarrow w^{(1)}_2 \rightarrow \cdots \rightarrow
w^{(1)}_{N-1}\rightarrow w^{(2)}_1\rightarrow\cdots
\rightarrow w^{(K_{N-3}+1)}_{N-1} \rightarrow 
\cdots \rightarrow w^{(K_{N-2})}_{N-1},
$$
that is, first in $w^{(1)}_1$, next in $w^{(1)}_2$ etc.
After the integration in $w$ 
we integrate in the variable $v$ by the order 
$$
v^{(1)}_1\rightarrow v^{(1)}_2 \rightarrow 
\cdots \rightarrow v^{(1)}_{N-1}
\rightarrow v^{(2)}_1 \rightarrow \cdots
\rightarrow v^{(K_{N-3}+1)}_{N-1} \rightarrow 
\cdots \rightarrow v^{(K_{N-2})}_{N-1}.
$$

In the variable $w^{(1)}_1$ the poles of the differential
form in the integrand outside
the contour $C^{(1)}_1$ is only at $w^{(1)}_1=q^{-1}w^{(1)}_2$.
It means that there are no poles at infinity too.
Thus we can calculate the integral in $w^{(1)}_1$ by
taking the residue at $w^{(1)}_1=q^{-1}w^{(1)}_2$.
After taking this residue the integrand have the same
structure in the variable $w^{(1)}_2$ and so on.
Therefore the integral in $w$`s is calculated by taking 
residues successively.
After calculating the integral in the variables $w$
the poles of the integrand in the variable $v^{(1)}_1$
inside the contour $\tilde{C}^{(1)}_1$ is only at 
$v^{(1)}_1=q^{-1}v^{(1)}_2$.
Hence the integral is calculated by taking the residue
at $v^{(1)}_1=q^{-1}v^{(1)}_2$.
After taking the residue in $v^{(1)}_1$ 
the integrand has the same structure
in the variable $v^{(1)}_2$ and so on.

Thus the integral is calculated by substituting
\begin{eqnarray}
&&
{
(q^{-1}-q)w^{(a)}_k \over 
(w^{(a)}_k-qw^{(a)}_{k+1})
(w^{(a)}_k-q^{-1}w^{(a)}_{k+1})
}=q^{-1},
\quad
w^{(a)}_k=q^{k+1}z_a,
\nonumber
\\
&&
{
(q^{-1}-q)v^{(b)}_{k+1} \over 
(v^{(b)}_k-qv^{(b)}_{k+1})
(v^{(b)}_k-q^{-1}v^{(b)}_{k+1})
}=1,
\quad
v^{(b)}_k=q^{k+1}u_b,
\nonumber
\end{eqnarray}
into the integrand and multiplying it by
$$
(-1)^{\sum_{r=0}^{N-1}(N-1-r)k_r}
$$
which comes from taking the residue in $w$ outside the
contour.
\vskip5mm

\noindent
{\bf Example.} $m=n=1$ case.
\noindent

In this case $i=j$, $\ep_1=\mu_1$, $k_r=l_r$ for any $r$ and
$r_0=0$. We consider the case $i=0$.
The integral formula is read as
\begin{eqnarray}
&&
\bar{G}^{(00)}(\zeta_1 \vert \xi_1)^\ep_\ep
=
C^{(00)}(\ep \vert \ep)\zeta_1^{-\ep}\xi_1^{\ep}\times I,
\nonumber
\\
&&
I=
\int_{
C^{(1)}_{\ep+1}
}
{d w^{(1)}_{\ep+1} \over 2\pi i}
\cdots
\int_{
C^{(1)}_{N-1}
}
{d w^{(1)}_{N-1} \over 2\pi i}
\int_{
\tilde{C}^{(1)}_{\ep+1}
}
{d v^{(1)}_{\ep+1} \over 2\pi i}
\cdots
\int_{
\tilde{C}^{(1)}_{N-1}
}
{d v^{(1)}_{N-1} \over 2\pi i}
\nonumber
\\
&&
\prod_{k=\ep+1}^{N-1}
{
(q^{-1}-q) w^{(1)}_k
\over
(w^{(1)}_k-q^{-1}w^{(1)}_{k+1})
(w^{(1)}_k-qw^{(1)}_{k+1})
}
\prod_{k=\ep+1}^{N-1}
{
(q^{-1}-q)v^{(1)}_{k+1}
\over
(v^{(1)}_k-q^{-1}v^{(1)}_{k+1})
(v^{(1)}_k-qv^{(1)}_{k+1})
}
\nonumber
\\
&&
\times
\prod_{k=\ep+1}^{N-1}(w^{(1)}_k-v^{(1)}_{k+1})
\prod_{k=\ep+2}^{N}(v^{(1)}_{k-1}-w^{(1)}_{k})
\prod_{k=\ep+1}^{N-1}
{1 \over
(w^{(1)}_k-qv^{(1)}_{k})
(w^{(1)}_k-q^{-1}v^{(1)}_{k})
}.
\nonumber
\end{eqnarray}
Here, by calculation, $C^{(00)}(\ep \vert \ep)=1$.
Let us denote the integrand of $I$ by $J$.
Consider the integral in $w^{(1)}_1$.
By definition of the contour $C^{(1)}_{\ep+1}$,
$q^{-1}w^{(1)}_{\ep+2}$ are outside and all other
poles on the complex plane are inside of $C^{(1)}_{\ep+1}$.
The differential form $Jd w^{(1)}_{\ep+1}$ has no poles at $\infty$.
We have
\begin{eqnarray}
&&
\hbox{Res}_{w^{(1)}_{\ep+1}=q^{-1}w^{(1)}_{\ep+2}}
J d w^{(1)}_{\ep+1}
=
\nonumber
\\
&&
-
\prod_{k=\ep+2}^{N-1}
{
(q^{-1}-q)w^{(1)}_k
\over
(w^{(1)}_k-q^{-1}w^{(1)}_{k+1})
(w^{(1)}_k-qw^{(1)}_{k+1})
}
\prod_{k=\ep+1}^{N-1}
{
(q^{-1}-q)v^{(1)}_{k+1}
\over
(v^{(1)}_k-q^{-1}v^{(1)}_{k+1})
(v^{(1)}_k-qv^{(1)}_{k+1})
}
\nonumber
\\
&&
\times
\prod_{k=\ep+2}^{N-1}(w^{(1)}_k-v^{(1)}_{k+1})
\prod_{k=\ep+3}^{N}(v^{(1)}_{k-1}-w^{(1)}_{k})
\prod_{k=\ep+3}^{N-1}
{1 \over
(w^{(1)}_k-qv^{(1)}_{k})
(w^{(1)}_k-q^{-1}v^{(1)}_{k})
}
\nonumber
\\
&&
{1
\over
(w^{(1)}_{\ep+2}-q^2v^{(1)}_{\ep+1})
(w^{(1)}_{\ep+2}-q^{-1}v^{(1)}_{\ep+2})
}.
\label{w1}
\end{eqnarray}

Consider this function (\ref{w1}) as a function of 
$w^{(1)}_{\ep+2}$.
By the definition of the contour $C^{(1)}_{\ep+2}$,
$q^{-1}w^{(1)}_{\ep+3}$ is outside of $C^{(1)}_{\ep+2}$
and all other poles on the complex plane are inside.
The differential forms (\ref{w1})$\times d w^{(1)}_{\ep+2}$
has no singularity at $\infty$.
Thus 
$$
\int_{C^{(1)}_{\ep+2}}\hbox{(\ref{w1})}d w^{(1)}_{\ep+2}
=-Res_{w^{(1)}_{\ep+2}=q^{-1}w^{(1)}_{\ep+3}}
\hbox{(\ref{w1})}d w^{(1)}_{\ep+2}
$$
and so on.
Consequently we have
\begin{eqnarray}
&&
\int_{
C^{(1)}_{\ep+1}
}
{d w^{(1)}_{\ep+1} \over 2\pi i}
\cdots
\int_{
C^{(1)}_{N-1}
}
{d w^{(1)}_{N-1} \over 2\pi i}
J
\nonumber
\\
&&
=(-1)^{N-1-\ep}
Res_{ w^{(1)}_{N-1}=q^N z_1 }
Res_{ w^{(1)}_{N-2}=q^{-1}w^{(1)}_{N-1} }
\cdots
Res_{ w^{(1)}_{\ep+1}=q^{-1}w^{(1)}_{\ep+2} }
J d w^{(1)}_{\ep+1}\cdots d w^{(1)}_{N-1}
\nonumber
\\
&&
=
q^{\ep+1}
{
z_1-qu_i
\over
q^{\ep+1}z_1-v^{(1)}_{\ep+1}
}
\prod_{k=\ep+1}^{N-1}
{
(q^{-1}-q)v^{(1)}_{k+1}
\over
(v^{(1)}_k-q^{-1}v^{(1)}_{k+1})
(v^{(1)}_k-qv^{(1)}_{k+1})
}.
\label{Jw}
\end{eqnarray}
A similar consideration is applicable to the function
(\ref{Jw}) in the variables $v$`s.
Finally we have
\begin{eqnarray}
I&=&
Res_{ v^{(1)}_{N-1}=q^{N}u_1 }
\cdots
Res_{ v^{(1)}_{\ep+1}=q^{-1}v^{(1)}_{\ep+2} }
\hbox{(\ref{Jw})}d v^{(1)}_{\ep+1}\cdots d v^{(1)}_{N-1}
\nonumber
\\
&=&1.
\nonumber
\end{eqnarray}
Consequently 
\begin{eqnarray}
&&
\bar{G}^{(00)}(\zeta_1 \vert \xi_1)^\ep_\ep
=(\zeta_1^{-1}\xi_1)^\ep.
\nonumber
\end{eqnarray}
This reproduces the formula in \cite{DO}.

%%%%%%%%%%%%%%%%%%%%%%%%%%%%%%%%%%%%%%%%%%%%%%%%%%%%%%%%%%%%%%%%%
\section{Integral formula for the trace of 
intertwining operators}\label{trf}
We recall the definition of $\bar{G}^{(i)}$: 
$$
\bar{G}^{(i)}({\bf \zeta} \vert {\bf \xi}\vert x,y)
^{\ep_1,\cdots,\ep_m}
_{\mu_1,\cdots,\mu_n}
=
{
\hbox{tr}_{V(\La_i)}\big(
x^Dy^H
\Phi_{\ep_1}(\zeta_1)
\cdots
\Phi_{\ep_{m}}(\zeta_m)
\Psi^{\ast}_{\mu_n}(\xi_n)
\cdots
\Psi^{\ast}_{\mu_1}(\xi_1)
\big)
\over
\bar{F}(z \vert u\vert x)
},
$$
where $\bar{F}(z \vert u\vert x)$ is given by (\ref{fun1}) 
and (\ref{fun2}).

We define
\begin{eqnarray}
&&
\bar{A}_r=\{j \vert \ep_j=r\},
\quad
A_r=\bar{A}_0\sqcup\cdots \sqcup \bar{A}_r,
\quad
\bar{B}_r=\{j \vert \mu_j=r\},
\quad
B_r=\bar{B}_0\sqcup\cdots \sqcup \bar{B}_r.
\nonumber
\end{eqnarray}
Then
$\sharp \bar{A}_r=k_r$, $\sharp \bar{B}_r=l_r$
and $\sharp A_r=K_r$, $\sharp B_r=L_r$.

The condition that the weight, with respect to $\uqpn$,
of the composition of the intertwining operators are zero is
$$
k_r-l_r=k_{N-1}-l_{N-1}=:r_0
$$
for $0\leq r\leq N-1$. We assume this condition.
We set $(z)_\infty=(z;x^N)_\infty$, 
$z_a=\zeta_a^N$, $u_b=\xi_b^N$. 
Then for $0\leq i\leq N-1$ we have

$$
\bar{G}^{(i)}({\bf \zeta} \vert {\bf \xi})
^{\ep_1,\cdots,\ep_m}
_{\mu_1,\cdots,\mu_n}
\qquad\qquad\qquad\qquad\qquad\qquad\qquad\qquad\qquad\qquad
\qquad\qquad\qquad\qquad\qquad\qquad\qquad
$$

$$
=
C^{tr(i)}(\ep\vert\mu)
\prod_{a=1}^m
\zeta_a^{(N-1)(m-n-a+1)-\ep_a+i}
\prod_{b=1}^n
\xi_b^{(N-1)(b-1)+\mu_b-i}
\prod_{a\in A_{N-2}}z_a^{-1}
\qquad\qquad\qquad\qquad\qquad\qquad
$$

$$
\times
\prod_{a<b,b\in A_{N-2}}z_a^{-1}
\prod_{a<b,a\in B_{N-2}}u_b^{-1}
\prod_{a\in A_{N-2},k}
\int_{C^{tr(a)}_k} {dw^{(a)}_k \over 2\pi i w^{(a)}_k}
\prod_{b\in B_{N-2},k}
\int_{\tilde{C}^{tr(b)}_k} {dv^{(b)}_k \over 2\pi i v^{(b)}_k}
\qquad\qquad\qquad\qquad
$$

$$
\times
\prod_{a<b,a,b\in A_{N-2}}
(w^{(a)}_{\ep_b+1})^{\theta(\ep_a\leq \ep_b)}
(w^{(a)}_{\ep_b})^{-\theta(\ep_a< \ep_b)}
\prod_{a<b,a,b\in B_{N-2}}
(v^{(b)}_{\mu_a+1})^{\theta(\mu_a\geq \mu_b)}
(v^{(b)}_{\mu_a})^{-\theta(\mu_a> \mu_b)}
\qquad\qquad\quad\,
$$

$$
\times
\prod_{a\in A_{N-2},b\in B_{N-2}}
(w^{(a)}_{\mu_b})^{\theta(\ep_a<\mu_b)}
(v^{(b)}_{\ep_a})^{\theta(\ep_a>\mu_b)}
\prod_{a\in A_{N-2}}
(w^{(a)}_{N-1})^{l_{N-1}}
\prod_{b\in B_{N-2}}
(v^{(b)}_{N-1})^{k_{N-1}}
\qquad\qquad\qquad\,\,\,
$$

$$
\times
\prod_{a<b,a\in A_{N-2},b\in \bar{A}_{N-1}}(w^{(a)}_{N-1})^{-1}
\prod_{a<b,b\in B_{N-2},a\in \bar{B}_{N-1}}(v^{(b)}_{N-1})^{-1}
\prod_{a\in A_{N-2}}w^{(a)}_{\ep_a+1}
\qquad\qquad\qquad\qquad\qquad
$$

$$
\times
\prod_{a<b}\prod_{k}(1-qw^{(a)}_k/w^{(b)}_{k+1})
\prod_{a>b}\prod_{k}(1-qw^{(b)}_{k+1}/w^{(a)}_k)
\prod_{a,b,k}
{1\over
(qw^{(a)}_k/w^{(b)}_{k+1})_\infty
(qw^{(b)}_{k+1}/w^{(a)}_{k})_\infty
}
\qquad\quad\,\,
$$

$$
\times
\prod_{a<b}\prod_{k}(1-q^{-1}v^{(b)}_{k+1}/v^{(a)}_k)
\prod_{a>b}\prod_{k}(1-q^{-1}v^{(a)}_k/v^{(b)}_{k+1})
\prod_{a,b,k}
{1\over
(q^{-1}v^{(a)}_k/v^{(b)}_{k+1})_\infty
(q^{-1}v^{(b)}_{k+1}/v^{(a)}_{k})_\infty
}
\quad\,\,
$$

$$
\times
\prod_{a,b,k}
{
\theta_{x^N}(v^{(b)}_{k+1}/w^{(a)}_k)
\over
(x^N)_\infty
}
\prod_{a,b,k}
{
\theta_{x^N}(w^{(a)}_{k+1}/v^{(b)}_k)
\over
(x^N)_\infty
}
\prod_{a,b}\prod_{k\leq N-1}
{
(x^N)_\infty^2
\over
\theta_{x^N}(qv^{(b)}_k/w^{(a)}_k)
\theta_{x^N}(qw^{(a)}_k/v^{(b)}_k)
}
\qquad\qquad
$$

$$
\times
\prod_{a<b}\prod_{k\leq N-1}
{
\theta_{x^N}(w^{(b)}_k/w^{(a)}_k)
\over
(x^N)_\infty
}
(q^2x^Nw^{(a)}_k/w^{(b)}_k)_\infty
(q^2w^{(b)}_k/w^{(a)}_k)_\infty
\qquad\qquad\qquad\qquad\qquad\qquad\qquad\!\!
$$

$$
\times
\prod_{a<b}\prod_{k\leq N-1}
{
\theta_{x^N}(v^{(a)}_k/v^{(b)}_k)
\over
(x^N)_\infty
}
(q^{-2}v^{(a)}_k/v^{(b)}_k)_\infty
(q^{-2}x^N v^{(b)}_k/v^{(a)}_k)_\infty
\qquad\qquad\qquad\qquad\qquad\qquad\qquad\!\!
$$

$$
\times
\big(
y_i \prod_{a:\ep_a+1\leq i}(w^{(a)}_i)^{-1}
\prod_{b:\mu_b+1\leq i}v^{(b)}_i 
\big)^{1-\delta_{i0}}
\theta_i(
g_0^{-1}g_1^2g_2^{-1},
\cdots,
g_{N-2}^{-1}g_{N-1}^2g_N^{-1}\vert x^N),
\qquad\qquad\qquad\qquad\!\!
$$
where
\begin{eqnarray}
&&
g_0^{-1}=(-1)^{(m-n)(N-1)},
\quad
g_N^{-1}=q^{(m-n)(N+1)}
\prod_{a=1}^m z_a\prod_{b=1}^n u_b^{-1},
\nonumber
\\
&&
g_j=
y_j\prod_{a:\ep_a+1\leq j}(w^{(a)}_j)^{-1}
\prod_{b:\mu_b+1\leq j}v^{(b)}_j
\quad
(1\leq j\leq N-1),
\nonumber
\\
&&
\theta_i(z_1,\cdots,z_{N-1}\vert p)
=
\sum_{\alpha\in \bar{Q}}
p^{{1\over2}(\alpha\vert \alpha)+(\alpha\vert \La_i)}
\prod_{j=1}^{N-1}z_j^{(\alpha\vert \La_j)}.
\nonumber
\end{eqnarray}
Here $\bar{Q}={\bf Z}\alpha_1\oplus\cdots\oplus{\bf Z}\alpha_{N-1}$
is the root lattice of $sl_N$.
The constant $C^{tr(i)}(\ep\vert \mu)$ is given by
\begin{eqnarray}
&&
C^{tr(i)}(\ep\vert \mu)
\nonumber
\\
&&
=(-1)^{
ir_0(N-1)+\hbox{sgn}_N(i)
+{1\over3}r_0(N-1)(N-2)(N-3)
+{1\over2}r_0(r_0+1)(N-1)
+K_{N-2}+nL_{N-2}
}
\nonumber
\\
&&
\times
(-1)^{
\sum_{a=1}^{N-2}ak_a
+\sum_{a\in A_{N-2}}a
+\sum_{b\in B_{N-2}}b
+\sum_{a\in A_{N-2},b\in B_{N-2}}c_{ab}
}
\nonumber
\\
&&
\times
(-1)^{
\sum_{a<b,a,b\in A_{N-2}}
(\ep_{ab}+\theta(\ep_a\leq \ep_b))
+
\sum_{a<b,a,b\in B_{N-2}}
(\mu_{ab}+\theta(\mu_a\geq \mu_b))
}
\nonumber
\\
&&
\times
q^{
ir_0(N+1)
+{1\over2}r_0^2N(N-1)(N+1)
-n(N+1)L_{N-2}
+(N-1)K_{N-2}L_{N-2}
+\sum_{b=1}^{N-2}bl_b
}
\nonumber
\\
&&
\times
q^{
(N+1)(-\sum_{a\in A_{N-2}}a+\sum_{b\in B_{N-2}}b)
-\sum_{a\in A_{N-2},b\in B_{N-2}}c_{ab}
}
\nonumber
\\
&&
\times
\Big(
{ \{q^2x^N\} \over \{q^{2N}x^N\} }
\Big)^m
\Big(
{ \{x^N\} \over \{q^{2N-2}x^N\} }
\Big)^n
(x^N)_\infty^{
\sum_{a=0}^{N-2}(N-1-a)(k_a+l_a)-1
}
\nonumber
\\
&&
\times
(q^2)_\infty
^{\sum_{a=0}^{N-2}(N-1-a)k_a}
(q^{-2})_\infty
^{\sum_{b=0}^{N-2}(N-1-b)l_b},
\nonumber
\end{eqnarray}
where we set
\begin{eqnarray}
&&
\ep_{ab}=\hbox{max}(\ep_a,\ep_b),
\quad
\mu_{ab}=\hbox{max}(\mu_a,\mu_b),
\quad
c_{ab}=\hbox{max}(\ep_a,\mu_b),
\nonumber
\\
&&
\hbox{sgn}_N(i)=
\nonumber
\\
&&
{1\over2}
\sum_{a=1}^m(N-N\{{i+a-1\over N}\})(N-1-N\{{i+a-1\over N}\})
\nonumber
\\
&&
+
{1\over2}
\sum_{b=1}^n(N-N\{{i+b-1\over N}\})(N-1-N\{{i+b-1\over N}\}).
\nonumber
\end{eqnarray}

The integration contour $C^{tr(a)}_k$ for $w^{(a)}_k$
and $\tilde{C}^{tr(b)}_k$ for $w^{(b)}_k$ 
are specified in the following manner.

The contour $C^{tr(a)}_k$ is a simple closed curve
going round the origin in the anticlockwise direction such that

\noindent
$qx^{Nm}w^{(b)}_{k\pm1}$ ($m\geq0$, any $b$),
$x^{Nm}w^{(b)}_k$ ($m\geq1$, $b\neq a$),
$x^{Nm}v^{(b)}_{k\pm1}$ ($m\geq1$, any $b$),
$q^{-1}x^{Nm}v^{(b)}_k$ ($m\geq0$, any $b$),
$q^{N+1}x^{Nm}u_b$ ($m\geq1$, any $b$),
$q^{N+2}x^{Nm}z_b$ ($m\geq0$, any $b$)
are inside and

\noindent
$q^{-1}x^{-Nm}w^{(b)}_{k\pm1}$ ($m\geq0$, any $b$),
$x^{-Nm}w^{(b)}_k$ ($m\geq1$, $b\neq a$),
$x^{-Nm}v^{(b)}_{k\pm1}$ ($m\geq1$, any $b$),
$qx^{-Nm}v^{(b)}_k$ ($m\geq1$, any $b$),
$q^{N+1}x^{-Nm}u_b$ ($m\geq1$, any $b$),
$q^{N}x^{-Nm}z_b$ ($m\geq0$, any $b$)
are outside.

The contour $\tilde{C}^{tr(b)}_k$ 
is a simple closed curve
going round the origin in the anticlockwise direction such that

\noindent
$q^{-1}x^{Nm}v^{(a)}_{k\pm1}$ ($m\geq 0$, any $a$),
$q^{-2}x^{Nm}v^{(b)}_k$ ($m\geq 1$, $b\neq a$),
$q^{N+1}x^{Nm}z_a$ ($m\geq1$, any $a$),
$q^Nx^{Nm}u_a$ ($m\geq0$, any $a$),
$x^{Nm}w^{(a)}_{k\pm1}$ ($m\geq 1$, any $a$),
$q^{-1}x^{Nm}w^{(a)}_k$ ($m\geq 1$, any $a$)
are inside and

\noindent
$qx^{-Nm}v^{(a)}_{k\pm1}$ ($m\geq 0$, any $a$),
$q^{2}x^{-Nm}v^{(b)}_k$ ($m\geq 1$, $b\neq a$),
$q^{N+1}x^{-Nm}z_a$ ($m\geq1$, any $a$),
$q^{N+2}x^{-Nm}u_a$ ($m\geq1$, any $a$),
$x^{-Nm}w^{(a)}_{k\pm1}$ ($m\geq 1$, any $a$),
$qx^{-Nm}w^{(a)}_k$ ($m\geq 1$, any $a$)
are outside.

It can be checked that 
those contours are well defined for $|x|<|q|^{2/N}<1$.

If we set $x=0$ in the formula above, we obtain the integral 
formula for the matrix element with $i=j$ in section \ref{matint}.

We have verified that this formula coincides with
the trace formula in \cite{JM} for $N=2$.

The derivation of the integral formula
of the trace is totally similar to the $sl_2$
case \cite{JMMN}\cite{JM} and it is briefly
explained in the appendix.

As a corollary of the integral formula for the trace 
we have

\begin{cor}
The functions 
$G({\bf \zeta} \vert {\bf \xi} \vert x,{\bf y})$
and 
$\bar{G}({\bf \zeta} \vert {\bf \xi} \vert x,{\bf y})$
are meromorphic functions on $({\bf C}^{\ast})^{n+m}$, where
${\bf C}^{\ast}={\bf C}\backslash \{0\}$
is the algebraic torus.
\end{cor}

\noindent
{\it Proof}. The singularity of the integral appears
only when the pinch of the integration contour occurs.
By the definition of the contour the pinch happens 
at $X=q^ax^bY$ for some integers $a,b$, 
where $X, Y\in\{\zeta_1,\cdots,\zeta_m,\xi_1,\cdots,\xi_n\}$.
Suppose that pinch occurs at $X=q^ax^bY$.
We decompose  the integral into the sum of residues and
the integral with the integration contour for which 
the pinch does not occur at $X=q^ax^bY$.
Since the integrand of the trace formula
is a meromorphic function on $({\bf C}^{\ast})^{n+m}$
its residue is also a meromorphic function
on $({\bf C}^{\ast})^{n+m}$.
In the decomposition the singularity at $X=q^ax^bY$ appears
only from the residue part.
Thus the singularity of the trace function at $X=q^ax^bY$
is a pole.

%%%%%%%%%%%%%%%%%%%%%%%%%%%%%%%%%%%%%%%%%%%%%%%%%%%%%%%%%%%%%%%%%%
\section{Discussion}
In this paper we have proved that the
trace of the composition of the intertwining operators of type I and
type II gives a basis of the solution space of the
qKZ equation at generic values of parameters.
The qKZ equation considered in this paper takes the value
in the tensor product of the finite dimensional 
irreducible $\uqnf$ module with the highest weight
$\La_1$.

There is a problem whether it is possible to construct solutions
of the qKZ equation taking values in the tensor product
of the arbitrary finite dimensional irreducible 
$\uqnf$ modules as a trace of intertwining
operators.
For $N=2$ it will be possible to construct solutions 
of the qKZ equation taking values in the arbitrary
irreducible $U_q(sl_2)$ modules by taking the
trace of the intertwining operators introduced in \cite{Na}
over the tensor product of integrable highest weight
$U_q(\widehat{sl_2})$ modules of level one \cite{HKMW}.
It is natural to expect that the trace functions
thus constructed 
give a basis of the solution space.
For $N\geq 3$ a similar construction will be possible.
For the moment what kind of modules we can treat is not very clear.
%%%%%%%%%%%%%%

Let us consider the qKZ equation (\ref{hqkz}) of $N=2$
on the weight subspace of 
$V_1 \ot \cdots \ot V_n$ with a weight, say $\lambda$.
At some special values of $\kappa$, which depend on
$p$, $q$ and $\lambda$,
the hypergeometric solution of Tarasov-Varchenko \cite{TV}
takes the value in the space of singular vectors 
with respect to certain action of $U_q(sl_2)$.
From the experience of rational and level zero case \cite{NPT},
it is probable that the trace function still gives 
a basis of the full space of the tensor product
at those special values of $\kappa$.
This means that the hypergeometric solution 
and the trace solution have
very different structures.
It is an interesting and important problem
to relate these two basis.
A partial result in this direction
is given in \cite{NPT}.

One of the important properties of the trace construction
of the solution is that it gives a map from 
$V^{\ot n}$ with fixed $n$
to the space of solutions of the qKZ equation taking the value in
$V^{\ot m}$ for any $m$.
This will be a key structure to relate finite and infinite
dimensional modules.
Note that it is nothing but the typical structure
of the form factors in integrable quantum field theories
\cite{S}\cite{NPT}.
The above mentioned problem connecting two types of solutions
is also important to understand the completeness problem
of local fields constructed by Smirnov \cite{S}\cite{BBS}.

The value $x=q^2$ is of particular interest,
since the correlation functions
and the form factors of the solvable lattice model
are given by some special case of the trace function
at this value of $x$.
We conjecture that $\det \bar{G}$ does not vanish identically
at $x=q^2$.

The generalization of the results in this paper
to other types of quantum affine algebra
is also interesting.

The trace of intertwining operators are also
studied in \cite{E}. 
Here we simply comments the following things.
In \cite{E} the trace is twisted by the Dynkin 
diagram automorphism and thus it is different from
the trace considered in this paper.
The difference equations satisfied by the trace in \cite{E} and
in this paper are also different.
\vskip1cm

\noindent
{\bf Acknowledgement}
I would like to thank Vitaly Tarasov for the helpful discussion.
This work is done while the author stays at LPTHE in
Universite Pierre et Marie Curie. I am grateful to
the laboratoire, in particular, to 
Olivier Babelon and Fedor Smirnov for their kind hospitality.

\appendix
\section
{Integral formula for trace
-$U_q(\widehat{\hbox{sl}_2})$ case-}
In the case of $U_q(\widehat{\hbox{sl}_2})$
the integral formula for the trace is given in \cite{JM}.
Our formula at $N=2$ in section \ref{trf} recovers it.
In this case the sum of the trace over $V(\La_0)$ and
$V(\La_1)$ simplifies a bit. It is used in the calculation of the
example in section 7. Thus we shall present this simplified formula.
It also serves as a simplest example of the trace formula.

\begin{eqnarray}
&&
{\bar G}({\bf \zeta}\vert{\bf \xi}\vert x,y)
^{\ep_1,\cdots,\ep_m}
_{\mu_1,\cdots,\mu_n}
=
\nonumber
\\
&&
C^{mn}_{st}
\prod_{j=1}^m\zeta_j^{-j+(1+\ep_j)/2}
\prod_{k=1}^n\xi_k^{k-(1+\mu_k)/2}
\prod_{r=1}^s \int_{{\cal C}} {dw_r \over 2\pi i w_r}
\prod_{r=1}^t \int_{\tilde{\cal C}} {dv_r \over 2\pi i v_r}
F^{AB}(\zeta,\xi,w,v\vert x,y),
\nonumber
\end{eqnarray}
where
\begin{eqnarray}
&&
F^{AB}(\zeta,\xi,w,v \vert x,y)=
\nonumber
\\
&&
\prod_{r=1}^s
\big(
\prod_{j<a_r}(qz_j-q^{-1}w_r)
\prod_{j>a_r}(z_j-w_r)
\big)
\prod_{j,r}{1\over (w_r/z_j)_\infty(q^2z_j/w_r)_\infty}
\nonumber
\\
&&
\times
\prod_{r=1}^t
\big(
\prod_{k<b_r}(qv_r^{-1}-q^{-1}u_k^{-1})
\prod_{k>b_r}(v_r^{-1}-u_k^{-1})
\big)
\prod_{k,r}{1\over (q^{-2}v_r/u_k)_\infty(u_k/v_r)_\infty}
\nonumber
\\
&&
\times
\prod_{r,j}{\theta_{x^2}(-q^{-1}v_r/z_j) \over (x^2)_\infty}
\prod_{r,k}{\theta_{x^2}(-qu_k/w_r) \over (x^2)_\infty}
\nonumber
\\
&&
\times
\prod_{r,r^\prime}
{
(x^2)_\infty^2
\over
\theta_{x^2}(-qv_r/w_{r^\prime})
\theta_{x^2}(-q^{-1}v_r/w_{r^\prime})
}
\nonumber
\\
&&
\times
\prod_{r<r^\prime}
{
(q^2w_r/w_{r^\prime})_\infty
(q^2w_{r^\prime}/w_r)_\infty
\over
w_{r^\prime}-q^2w_r
}
{
w_{r^\prime}^{-1}
\theta_{x^2}(w_{r^\prime}/w_r)
\over
(x^2)_\infty
}
\nonumber
\\
&&
\prod_{r<r^\prime}
(v_{r^\prime}-q^{-2}v_r)
(x^2q^{-2}v_r/v_{r^\prime})_\infty
(x^2q^{-2}v_{r^\prime}/v_r)_\infty
{
v_{r^\prime}
\theta_{x^2}(v_r/v_{r^\prime})
\over
(x^2)_\infty
}
\nonumber
\\
&&
\times
\theta_x(
(-1)^{t+1}
(-q)^{m-n\over2}y
{
\prod\zeta_j\prod v_r
\over
\prod\xi_k\prod w_r
}
).
\nonumber
\end{eqnarray}

Here 
$A=\{a_1<\cdots<a_s\}=\{j \vert \ep_j=+\}$ and
$B=\{b_1<\cdots<b_t\}=\{j \vert \mu_j=+\}$,
$z_j=\zeta_j^2$, $u_k=-\xi_k^2$.
The integral contour ${\cal C}$ and $\tilde{\cal C}$ 
go round the origin such that

\noindent
for ${\cal C}$: $q^2x^{2l}z_j$ $(l\geq 0)$ are inside and 
$x^{-2l}z_j$ $(l\geq 0)$ are outside,

\noindent
for $\tilde{\cal C}$: $x^{2l}u_k$ $(l\geq 0)$ are inside and 
$q^2x^{-2l}u_k$ $(l\geq 0)$ are outside,
$-q^{\pm1}x^{2l}w_r$ $(l\geq 1)$ are inside and 
$-q^{\pm1}x^{-2l}w_r$ $(l\geq0)$ are outside.

We have rewritten the formula in \cite{JM} using the following formula:
\begin{eqnarray}
\theta_{x^4}(-xX^2)+(-1)^tX\theta_{x^4}(-x^3X^2)
=
\theta_x((-1)^{t+1}X).
\nonumber
\end{eqnarray}

\section{Boson expression of intertwining operators}
Here we recall the bosonic expression of
intertwining operators for $\uqn$ \cite{Ko}.

Let us consider the Heisenberg algebra generated by
$\{a_i(k) \vert 1\leq i\leq N-1, k\in {\bf Z}\backslash\{0\}\}$
with the commutation relation
\begin{eqnarray}
&&
[a_i(k), a_j(l)]=
\delta_{k+l,0}
{[(\alpha_i \vert \alpha_j)k][k] \over k}.
\nonumber
\end{eqnarray}
Let ${\cal H}$ be the Fock space of this algebra,
${\cal H}={\bf C}[a_i(-k)\vert 1\leq i \leq N-1,
k\in {\bf Z}\backslash\{0\}]$.
Let $\bar{Q}=\oplus_{j=1}^{N-1}{\bf Z}\alpha_j$ be the root
lattice of $sl_N$.
Then the twisted group algebra ${\bf C}[\bar{Q}]$ is the algebra
generated by $\hbox{e}^{\alpha_1}$, ... , $\hbox{e}^{\alpha_{N-1}}$
with the defining relation
\begin{eqnarray}
&&
\hbox{e}^{\alpha_i}\hbox{e}^{\alpha_j}
=
(-1)^{(\alpha_i \vert \alpha_j)}
\hbox{e}^{\alpha_j}\hbox{e}^{\alpha_i}.
\nonumber
\end{eqnarray}

Then

\begin{theo}\cite{FJ}
There is an isomorphism 
\begin{eqnarray}
&&
V(\La_i)\simeq {\cal H}\otimes{\bf C}[\bar{Q}]\hbox{e}^{\La_i},
\label{fj}
\end{eqnarray}
where ${\cal H}\otimes{\bf C}[\bar{Q}]\hbox{e}^{\La_i}$
is the vector space consisting of the symbols 
$X\hbox{e}^{\La_i}$, $X\in {\bf C}[\bar{Q}]$.
\end{theo}

For the action of the generators of $\uqn$ on the right hand side
of (\ref{fj}) see \cite{FJ}\cite{Ko}.

To describe the intertwining operators we introduce
the algebra containing 
${\bf C}[\bar{Q}]\hbox{e}^{\La_i}$.
Let $\bar{P}=\oplus_{j=1}^{N-1}{\bf Z}\La_j=
\oplus_{j=2}^{N-1}{\bf Z}\alpha_j\oplus{\bf Z}\La_{N-1}$ be
the weight lattice of $sl_N$.
Then the extended group algebra ${\bf C}[\bar{P}]$ is the algebra
generated by 
$\hbox{e}^{\alpha_1}$, ... , $\hbox{e}^{\alpha_{N-1}}$,
$\hbox{e}^{\La_{N-1}}$ with the defining relation \cite{Ko}
\begin{eqnarray}
&&
\hbox{e}^{\alpha}\hbox{e}^{\beta}
=
(-1)^{(\alpha \vert \beta)}
\hbox{e}^{\beta}\hbox{e}^{\alpha},
\quad
\alpha,\beta\in
\{\alpha_2,\cdots,\alpha_{N-1},\La_{N-1}\}.
\nonumber
\end{eqnarray}
As a convention, 
for $\alpha=\sum_{j=2}^{N-1}m_j\alpha_j+m_N\La_{N-1}$,
we set
\begin{eqnarray}
&&
\hbox{e}^{\alpha}
=
\hbox{e}^{m_2\alpha_2}\cdots \hbox{e}^{m_{N-1}\alpha_{N-1}}
\hbox{e}^{m_{N}\La_{N-1}}.
\nonumber
\end{eqnarray}
Note that
\begin{eqnarray}
\alpha_1=-\sum_{r=2}^{N-1}r\alpha_r+N\La_{N-1},
\quad
\La_i=-\sum_{r=i+1}^{N-1}(r-i)\alpha_r+(N-i)\La_{N-1}.
\nonumber
\end{eqnarray}
The algebra ${\bf C}[\bar{Q}]$ becomes a subalgebra
of ${\bf C}[\bar{P}]$.
We consider ${\bf C}[\bar{Q}]\hbox{e}^{\La_i}$
as a subspace of ${\bf C}[\bar{P}]$.

We define the action of the symbols
$\partial_\alpha$, $\hbox{e}^{\alpha}$ $(\alpha\in \bar{Q})$,
and $d$ on the space 
${\cal H}\otimes{\bf C}[\bar{Q}]\hbox{e}^{\La_i}$.
Let $X=a_{j_1}(-n_1)\cdots a_{j_k}(-n_k)\in {\cal H}$,
$\hbox{e}^{\beta}
\in {\cal H}\otimes{\bf C}[\bar{Q}]\hbox{e}^{\La_i}$
and $Y=X\ot \hbox{e}^{\beta}$.
Then
\begin{eqnarray}
&&
\partial_\alpha Y=(\alpha\vert\beta)Y,
\quad
\hbox{e}^\alpha Y=X\otimes \hbox{e}^\alpha\hbox{e}^\beta,
\nonumber
\\
&&
dY=(-\sum_{r=1}^k n_k -{(\beta\vert\beta)\over2}
+{(\La_i\vert\La_i)\over2})Y.
\nonumber
\end{eqnarray}
Then the principal grading operator $D^{(i)}$ on
$V(\La_i)$ is given by
\begin{eqnarray}
&&
D^{(i)}=-\rho+{i(N-i) \over2},
\quad
\rho=Nd+{1\over2}\sum_{r=1}^{N-1}r(N-r)\partial_{\alpha_r}.
\nonumber
\end{eqnarray}

We set
\begin{eqnarray}
X^{\pm}_j(w)&=&
\exp\Big(
\pm\sum_{k=1}^\infty
{a_j(-k) \over [k]}q^{\mp{k\over2}}w^k
\Big)
\exp\Big(
\mp\sum_{k=1}^\infty
{a_j(k) \over [k]}q^{\mp{k\over2}}w^{-k}
\Big)
\hbox{e}^{\pm\alpha_j}w^{\pm\partial_{\alpha_j}},
\nonumber
\\
&=&\sum_{n\in {\bf Z}}x^{\pm}_{j,n}w^{-n-1},
\nonumber
\\
x^{\pm}_j&=&x^{\pm}_{j,0}.
\nonumber
\end{eqnarray}
Then

\begin{theo}\cite{Ko}
$$
\tilde{\Phi}^{h(i)}_{N-1}(z)
=
\exp\Big(
\sum_{k=1}^\infty
a_{N-1}^{\ast}(-k)q^{(N+{3\over2})k}z^k
\Big)
\exp\Big(
\sum_{k=1}^\infty
a_{N-1}^{\ast}(k)q^{-(N+{1\over2})k}z^{-k}
\Big)
$$
$$
\times
\hbox{e}^{\La_{N-1}}
(q^{N+1}z)^{\partial_{{\La}_{N-1}}+{N-1-i \over N}}
(-1)^{(N-1)(\partial_{{\La}_{1}}-{N-1-i \over N})}
(-1)^{{1\over2}(N-i)(N-1-i)},
$$
$$
\tilde{\Phi}^{h(i)}_{j}(z)=
[\tilde{\Phi}^{h(i)}_{j+1}(z),x^{-}_j]_q,
\quad
0\leq j\leq N-2,
$$
%%%%%%%%%%%%%%%%%%%
%$$
%X^{+}_j(v)=
%\exp\Big(
%\sum_{k=1}^\infty
%{a_j(-k) \over [k]}q^{-{k\over2}}v^k
%\Big)
%\exp\Big(
%-\sum_{k=1}^\infty
%{a_j(k) \over [k]}q^{-{k\over2}}v^{-k}
%\Big)
%\hbox{e}^{\alpha_j}v^{\partial_{\alpha_j}},
%$$
$$
\tilde{\Psi}^{\ast h(i)}_{N-1}(u)
=
\exp\Big(
-\sum_{k=1}^\infty
a_{N-1}^{\ast}(-k)q^{(N+{1\over2})k}u^k
\Big)
\exp\Big(
-\sum_{k=1}^\infty
a_{N-1}^{\ast}(k)q^{-(N+{3\over2})k}u^{-k}
\Big)
$$
$$
\times
\hbox{e}^{-{\La}_{N-1}}
(q^{N+1}u)^{-\partial_{{\La}_{N-1}}+{i\over N}}
(-1)^{(N-1)(-\partial_{{\La}_{1}}+{N-i\over N})}
(-1)^{{1\over2}(N-i)(N-1-i)},
$$
$$
\tilde{\Psi}^{\ast h(i)}_{N-1}(u)=
[x^{+}_j,\tilde{\Psi}^{\ast h(i)}_{j+1}(u)]_{q^{-1}},
\quad
0\leq j\leq N-2,
$$
where
$[X,Y]_q=XY-qYX$ and
\begin{eqnarray}
&&
a_{N-1}^{\ast}(k)=
{-1 \over [k][Nk]}
\sum_{r=1}^{N-1}
[rk]a_r(k).
\nonumber
\end{eqnarray}
\end{theo}

The elements $a_{N-1}^{\ast}(k)$ satisfy the relations
$$
[a_j(k),a_{N-1}^{\ast}(-l)]
=\delta_{k,l}
\delta_{j,N-1}{[k] \over k},
\quad
[a_{N-1}^{\ast}(k),a_{N-1}^{\ast}(-l)]
=-\delta_{k,l}
{q^k(1-q^{(2N-2)k}) \over k(1-q^{2Nk})}.
$$
The inner product is given explicitly by
$$
({\La}_i \vert {\La}_j)
={i(N-j) \over N}\quad (i\leq j),
\quad
(\alpha_i \vert {\La}_{j})=\delta_{i,j}.
$$
%%%%%%%%%%%%%%%%%%%

\section{List of normal ordering rules}
We define the normal ordered operator as an operator
of the form
\begin{eqnarray}
&&
\exp\big(
\sum_{j=1}^{N-1}\sum_{n=1}^\infty
A^{(j)}_na_J(-n)
\big)
\exp\big(
\sum_{j=1}^{N-1}\sum_{n=1}^\infty
B^{(j)}_na_J(n)
\big)
\nonumber
\\
&&
\times
\exp(\sum_{j=1}^{N-1}c_j \alpha_j)
\exp(\sum_{j=1}^{N-1}c_j^\prime{\partial_{\alpha_j}})
\Big)
\nonumber
\end{eqnarray}
Thus we define the normal order of the product of operators
as
\begin{eqnarray}
:a_i(k)a_j(l):&=& a_i(k)a_j(l) \quad \hbox{ if $k\leq l$}
\nonumber
\\
              &=& a_j(l)a_i(k) \quad \hbox{ if $k>l$},
\nonumber
\\
:\partial_\alpha a_i(k):&=&:a_i(k) \partial_\alpha:=
a_i(k) \partial_\alpha,
\nonumber
\\
:\hbox{e}^\alpha a_i(k):&=&:a_i(k) \hbox{e}^\alpha:=
a_i(k) \hbox{e}^\alpha,
\nonumber
\\
:\partial_\alpha \hbox{e}^\beta:&=&:\hbox{e}^\beta\partial_\alpha:=
\hbox{e}^\beta\partial_\alpha.
\nonumber
\end{eqnarray}

We shall give a list of expressions of operators
in terms of their normal ordered operators.

$$
X^{-}_{j_1}(w_1)X^{-}_{j_2}(w_2)
=
:X^{-}_{j_1}(w_1)X^{-}_{j_2}(w_2):
\quad \vert j_1-j_2\vert>1, j_1,j_2\neq 0,
$$

$$
X^{-}_{1}(w_1)X^{-}_{j}(w_2)
=
(-1)^{j+1}
:X^{-}_{1}(w_1)X^{-}_{j}(w_2):
\quad  j\geq 3,
$$

$$
X^{-}_{j}(w_1)X^{-}_{1}(w_2)
=
(-1)^{j-1}
:X^{-}_{j}(w_1)X^{-}_{1}(w_2):
\quad j\geq 3,
$$

$$
X^{-}_{j}(w_1)X^{-}_{j+1}(w_2)
=
{(-1)^{\delta_{j1}} \over w_1-qw_2}
:X^{-}_{j}(w_1)X^{-}_{j+1}(w_2):,
$$

$$
X^{-}_{j+1}(w_1)X^{-}_{j}(w_2)
=
{(-1)^{1-\delta_{j1}} \over w_1-qw_2}
:X^{-}_{j+1}(w_1)X^{-}_{j}(w_2):,
$$

$$
X^{-}_{j}(w_1)X^{-}_{j}(w_2)
=
(w_1-q^2w_2)(w_1-w_2)
:X^{-}_{j}(w_1)X^{-}_{j}(w_2):,
$$
%%%%%%%%%%%

$$
X^{-}_{j_1}(w)X^{+}_{j_2}(v)
=
:X^{-}_{j_1}(w)X^{+}_{j_2}(v):
\quad \vert j_1-j_2\vert>1, j_,j_2\neq 1,
$$

$$
X^{-}_{1}(w_1)X^{+}_{j}(w_2)
=
(-1)^{j+1}
:X^{-}_{1}(w_1)X^{+}_{j}(w_2):
\quad  j\geq 3,
$$

$$
X^{-}_{j}(w_1)X^{+}_{1}(w_2)
=
(-1)^{j-1}
:X^{-}_{j}(w_1)X^{+}_{1}(w_2):
\quad  j\geq 3,
$$

$$
X^{-}_{j}(w)X^{+}_{j+1}(v)
=
(-1)^{\delta_{j1}}
(w-v)
:X^{-}_{j}(w)X^{+}_{j+1}(v):,
$$

$$
X^{-}_{j+1}(w)X^{+}_{j}(v)
=
(-1)^{1-\delta_{j1}}(w-v)
:X^{-}_{j+1}(w)X^{+}_{j}(v):,
$$

$$
X^{-}_{j}(w)X^{+}_{j}(v)
=
{1 \over (w-qv)(w-q^{-1}v)}
:X^{-}_{j}(w)X^{+}_{j}(v):,
$$

%%%%%%%

$$
X^{+}_{j_1}(v)X^{-}_{j_2}(w)
=
:X^{+}_{j_1}(v)X^{-}_{j_2}(w):
\quad \vert j_1-j_2\vert>1, j_1,j_2 \neq 1,
$$

$$
X^{+}_{1}(w_1)X^{-}_{j}(w_2)
=
(-1)^{j+1}
:X^{+}_{1}(w_1)X^{-}_{j}(w_2):
\quad j\geq 3,
$$

$$
X^{+}_{j}(w_1)X^{-}_{1}(w_2)
=
(-1)^{j-1}
:X^{+}_{j}(w_1)X^{-}_{1}(w_2):
\quad j\geq 3,
$$

$$
X^{+}_{j}(v)X^{-}_{j+1}(w)
=
(-1)^{\delta_{j1}}(v-w)
:X^{+}_{j}(v)X^{-}_{j+1}(w):,
$$

$$
X^{+}_{j+1}(v)X^{-}_{j}(w)
=
(-1)^{1-\delta_{j1}}(v-w)
:X^{+}_{j+1}(v)X^{-}_{j}(w):,
$$

$$
X^{+}_{j}(v)X^{-}_{j}(w)
=
{1 \over (v-qw)(v-q^{-1}w)}
:X^{+}_{j}(v)X^{-}_{j}(w):,
$$

%%%%%%%%
$$
X^{+}_{j_1}(v_1)X^{+}_{j_2}(v_2)
=
:X^{+}_{j_1}(v_1)X^{+}_{j_2}(v_2):
\quad \vert j_1-j_2\vert>1, j_1,j_2\neq 1,
$$

$$
X^{+}_{1}(w_1)X^{+}_{j}(w_2)
=
(-1)^{j+1}
:X^{+}_{1}(w_1)X^{+}_{j}(w_2):
\quad j\geq 3,
$$

$$
X^{+}_{j}(w_1)X^{+}_{1}(w_2)
=
(-1)^{j-1}
:X^{+}_{j}(w_1)X^{+}_{1}(w_2):
\quad j\geq 3,
$$

$$
X^{+}_{j}(v_1)X^{+}_{j+1}(v_2)
=
{(-1)^{\delta_{j1}} \over v_1-q^{-1}v_2}
:X^{+}_{j}(v_1)X^{+}_{j+1}(v_2):,
$$

$$
X^{+}_{j+1}(v_1)X^{+}_{j}(v_2)
=
{(-1)^{1-\delta_{j1}} \over v_1-q^{-1}v_2}
:X^{+}_{j+1}(v_1)X^{+}_{j}(v_2):,
$$

$$
X^{+}_{j}(v_1)X^{+}_{j}(v_2)
=
(v_1-q^{-2}v_2)(v_1-v_2)
:X^{+}_{j}(v_1)X^{+}_{j}(v_2):,
$$

%%%%%
$$
\tilde{\Phi}^{h(i)}_{N-1}(z)X^{-}_j(w)
=
:\tilde{\Phi}^{h(i)}_{N-1}(z)X^{-}_j(w):
\quad j\neq N-1,
$$

$$
\tilde{\Phi}^{h(i)}_{N-1}(z)X^{-}_{N-1}(w)
=
{q^{-1} \over w-q^Nz}
:\tilde{\Phi}^{h(i)}_{N-1}(z)X^{-}_{N-1}(w):,
$$

%%%
$$
\tilde{\Phi}^{h(i)}_{N-1}(z)X^{+}_j(v)
=
:\tilde{\Phi}^{h(i)}_{N-1}(z)X^{+}_j(v):
\quad j\neq N-1,
$$

$$
\tilde{\Phi}^{h(i)}_{N-1}(z)X^{+}_{N-1}(v)
=
(v-q^{N+1}z)
:\tilde{\Phi}^{h(i)}_{N-1}(z)X^{+}_{N-1}(v):,
$$

%%%%%%%%

$$
\tilde{\Psi}^{\ast h(i)}_{N-1}(u)X^{+}_j(v)
=
:\tilde{\Psi}^{\ast h(i)}_{N-1}(u)X^{+}_j(v):
\quad j\neq N-1,
$$

$$
\tilde{\Psi}^{\ast h(i)}_{N-1}(u)X^{+}_{N-1}(v)
=
{q \over v-q^{N+2}u}
:\tilde{\Psi}^{\ast h(i)}_{N-1}(u)X^{+}_{N-1}(v):,
$$

%%%%%%%%
$$
X^{-}_j(w)\tilde{\Phi}^{h(i)}_{N-1}(z)
=
:X^{-}_j(w)\tilde{\Phi}^{h(i)}_{N-1}(z):,
\quad j\neq N-1,
$$

$$
X^{-}_{N-1}(w)\tilde{\Phi}^{h(i)}_{N-1}(z)
=
{1 \over w-q^{N+2}z}
:X^{-}_{N-1}(w)\tilde{\Phi}^{h(i)}_{N-1}(z):,
$$
%%%%%%%%%%%

$$
X^{+}_j(v)\tilde{\Phi}^{h(i)}_{N-1}(z)
=
:X^{+}_j(v)\tilde{\Phi}^{h(i)}_{N-1}(z):,
\quad j\neq N-1,
$$

$$
X^{+}_{N-1}(v)\tilde{\Phi}^{h(i)}_{N-1}(z)
=
(v-q^{N+1}z)
:X^{+}_{N-1}(v)\tilde{\Phi}^{h(i)}_{N-1}(z):,
$$

%%%%%%%%%%
$$
X^{+}_j(v)\tilde{\Psi}^{\ast h(i)}_{N-1}(u)
=
:X^{+}_j(v)\tilde{\Psi}^{\ast h(i)}_{N-1}(u):,
\quad j\neq N-1,
$$

$$
X^{+}_{N-1}(v)\tilde{\Psi}^{\ast h(i)}_{N-1}(u)
=
{1 \over v-q^{N}u}
:X^{+}_{N-1}(v)\tilde{\Psi}^{\ast h(i)}_{N-1}(u):,
$$
%%%%%%%%%%%

$$
X^{-}_j(w)\tilde{\Psi}^{\ast h(i)}_{N-1}(u)
=
:X^{-}_j(w)\tilde{\Psi}^{\ast h(i)}_{N-1}(u):,
\quad j\neq N-1,
$$

$$
X^{-}_{N-1}(w)\tilde{\Psi}^{\ast h(i)}_{N-1}(u)
=
(w-q^{N+1}u)
:X^{-}_{N-1}(w)\tilde{\Psi}^{\ast h(i)}_{N-1}(u):,
$$

%%%%%%%%%%
$$
\tilde{\Phi}^{h(i_1)}_{N-1}(z_1)\tilde{\Phi}^{h(i_2)}_{N-1}(z_2)
=
(-q^{N+1}z_1)^{{N-1 \over N}}
{
(q^{2}{z_2 \over z_1})_\infty
\over
(q^{2N}{z_2 \over z_1})_\infty
}
:\tilde{\Phi}^{h(i_1)}_{N-1}(z_1)\tilde{\Phi}^{h(i_2)}_{N-1}(z_2):,
$$

$$
\tilde{\Phi}^{h(i_1)}_{N-1}(z)\tilde{\Psi}^{\ast h(i_2)}_{N-1}(u)
=
(-q^{N+1}z)^{-{N-1 \over N}}
{
(q^{2N-1}{u \over z})_\infty
\over
(q{u \over z})_\infty
}
:\tilde{\Phi}^{h(i_1)}_{N-1}(z)\tilde{\Psi}^{\ast h(i_2)}_{N-1}(u):,
$$

$$
\tilde{\Psi}^{\ast h(i_1)}_{N-1}(u)
\tilde{\Phi}^{h(i_2)}_{N-1}(z)
=
(-q^{N+1}u)^{-{N-1 \over N}}
{
(q^{2N-1}{z \over u})_\infty
\over
(q{z \over u})_\infty
}
:\tilde{\Psi}^{\ast h(i_1)}_{N-1}(u)
\tilde{\Phi}^{h(i_2)}_{N-1}(z):,
$$

$$
\tilde{\Psi}^{\ast h(i_2)}_{N-1}(u_2)
\tilde{\Psi}^{\ast h(i_1)}_{N-1}(u_1)
=
(-q^{N+1}u_2)^{{N-1 \over N}}
{
({u_1 \over u_2})_\infty
\over
(q^{2N-2}{u_1 \over u_2})_\infty
}
:\tilde{\Psi}^{\ast h(i_2)}_{N-1}(u_2)
\tilde{\Psi}^{\ast h(i_1)}_{N-1}(u_1):.
$$

Let us set
\begin{eqnarray}
&&
\tilde{\Phi}^{h(i)}_j(z \vert w_{N-1}\cdots w_{j+1})=
[\tilde{\Phi}^{h(i)}_{j+1}(z \vert w_{N-1}\cdots w_{j+2}), 
X^{-}_{j+1}(w_{j+1})]_q,
\nonumber
\\
&&
\tilde{\Psi}^{\ast h(i)}_j(u \vert v_{N-1}\cdots v_{j+1})=
[X^{+}_{j+1}(v_{j+1}),
\tilde{\Psi}^{\ast h(i)}_{j+1}(u \vert v_{N-1}\cdots v_{j+2})
]_{q^{-1}}.
\nonumber
\end{eqnarray}

Then
\begin{eqnarray}
&&
\tilde{\Phi}^{h(i)}_j(z \vert w_{N-1}\cdots w_{j+1})
\nonumber
\\
&&
=
(-1)^{{1\over2}N(N+1)\delta_{j0}+\delta_{j0}}
%{(q^{-1}-q)w_{N-1} \over (w_{N-1}-q^Nz)(w_{N-1}-q^{N+2}z)}
\prod_{k=j+1}^{N-1}
{(q^{-1}-q)w_k \over (w_{k}-q^{-1}w_{k+1})(w_k-qw_{k+1})}
\nonumber
\\
&&
\times
:\tilde{\Phi}^{h(i)}_{N-1}(z)
X^{-}_{N-1}(w_{N-1})
\cdots
X^{-}_{j+1}(w_{j+1}):,
\nonumber
\end{eqnarray}

\begin{eqnarray}
&&
\tilde{\Psi}^{\ast h(i)}_j(u \vert v_{N-1}\cdots v_{j+1})
\nonumber
\\
&&
=
(-1)^{{1\over2}N(N+1)\delta_{j0}+\delta_{j0}}
%{q^Nu(1-q^2) \over (v_{N-1}-q^Nu)(v_{N-1}-q^{N+2}u)}
\prod_{k=j+1}^{N-1}
{(q^{-1}-q)v_{k+1} \over (v_k-q^{-1}v_{k+1})(v_k-qv_{k+1})}
\nonumber
\\
&&
\times
:\tilde{\Psi}^{\ast h(i)}_{N-1}(u)
X^{+}_{N-1}(v_{N-1})
\cdots
X^{+}_{j+1}(v_{j+1}):,
\nonumber
\end{eqnarray}
where we set $w_N=q^{N+1}z$ and $v_N=q^{N+1}u$.

We have
$$
:\tilde{\Phi}^{h(i_1)}_{N-1}(z_1)
X^{-}_{N-1}(w^{(1)}_{N-1})
\cdots
X^{-}_{j_1+1}(w^{(1)}_{j_1+1}):
:\tilde{\Phi}^{h(i_2)}_{N-1}(z_2)
X^{-}_{N-1}(w^{(2)}_{N-1})
\cdots
X^{-}_{j_2+1}(w^{(2)}_{j_2+1}):
$$

$$
=
(-1)^{
{1\over2}(N-1-j_2)(N+2+j_2)\delta_{j_10}
+{1\over2}(N-1-j_1)(N+2+j_1)\delta_{j_20}
}
\qquad\qquad\qquad\qquad\qquad
$$

$$
\times
h^{(+)}({z_2 \over z_1})
q^{-1}
(-q^{N+1}z_1)^{{N-1 \over N}}
{
\prod_{k}
(w^{(1)}_k-q^2w^{(2)}_k)(w^{(1)}_k-w^{(2)}_k)
\over
(w^{(1)}_{N-1}-q^{N+2}z_2)(w^{(2)}_{N-1}-q^Nz_1)
}
\qquad\qquad
$$

$$
\times
\prod_{k}
{-1 \over w^{(1)}_k-qw^{(2)}_{k-1}}
\prod_{k}
{1 \over w^{(1)}_k-qw^{(2)}_{k+1}}
:\tilde{\Phi}^{h(i_1)}_{N-1}(z_1)
\cdots
X^{-}_{j_2+1}(w^{(2)}_{j_2+1}):,
$$
%%%%%%%%%%%%%%%%%%

$$
:\tilde{\Phi}^{h(i_1)}_{N-1}(z)
X^{-}_{N-1}(w_{N-1})
\cdots
X^{-}_{j_1+1}(w_{j_1+1}):
:\tilde{\Psi}^{\ast h(i_2)}_{N-1}(u)
X^{+}_{N-1}(v_{N-1})
\cdots
X^{+}_{j_2+1}(v_{j_2+1}):
$$

$$
=
(-1)^{
{1\over2}(N-1-j_2)(N+2+j_2)\delta_{j_10}
+{1\over2}(N-1-j_1)(N+2+j_1)\delta_{j_20}
}
\qquad\qquad\qquad\qquad\qquad
$$

$$
\times
h^{(0)}({u \over z})^{-1}
(-q^{N+1}z)^{-{N-1 \over N}}
{
(w_{N-1}-q^{N+1}u)(v_{N-1}-q^{N+1}z)
\over 
\prod_{k}(w_k-qv_k)(w_k-q^{-1}v_k)
}
\qquad\qquad
$$

$$
\times
\prod_{k}
(w_k-v_{k+1})
\prod_{k}
(-1)(w_k-v_{k-1})
:\tilde{\Phi}^{h(i_1)}_{N-1}(z)
\cdots
X^{+}_{j_2+1}(v_{j_2+1}):,
$$
%%%%%%%%%%%%%%%%%%%%%

$$
:\tilde{\Psi}^{\ast h(i_2)}_{N-1}(u_2)
X^{+}_{N-1}(v^{(2)}_{N-1})
\cdots
X^{+}_{j_2+1}(v^{(2)}_{j_2+1}):
:\tilde{\Psi}^{\ast h(i_1)}_{N-1}(u_1)
X^{+}_{N-1}(v^{(1)}_{N-1})
\cdots
X^{+}_{j_1+1}(v^{(1)}_{j_1+1}):
$$

$$
=
(-1)^{
{1\over2}(N-1-j_2)(N+2+j_2)\delta_{j_10}
+{1\over2}(N-1-j_1)(N-2+j_1)\delta_{j_20}
}
\qquad\qquad\qquad\qquad\qquad
$$

$$
\times
h^{(-)}({u_1 \over u_2})
q(-q^{N+1}u_2)^{{N-1 \over N}}
{
\prod_{k}
(v^{(2)}_k-q^{-2}v^{(1)}_k)
(v^{(2)}_k-v^{(1)}_k)
\over 
(v^{(1)}_{N-1}-q^{N+2}u_2)(v^{(2)}_{N-1}-q^Nu_1)
}
\qquad\qquad
$$

$$
\times
\prod_{k}
{1 \over v^{(2)}_{k}-q^{-1}v^{(1)}_{k+1}}
\prod_{k}
{-1 \over v^{(2)}_{k}-q^{-1}v^{(1)}_{k-1}}
:\tilde{\Psi}^{\ast h(i_2)}_{N-1}(u_2)
\cdots
X^{+}_{j_1+1}(v^{(1)}_{j_1+1}):.
$$

Here, denoting $(z)_\infty=(z;x^N)_\infty$, we set
$$
h^{(+)}(z)=
{(q^2z)_\infty \over (q^{2N}z)_\infty},
\quad
h^{(0)}(z)^{-1}=
{(q^{2N-1}z)_\infty \over (qz)_\infty},
\quad
h^{(-)}(z)=
{(z)_\infty \over (q^{2N-2}z)_\infty}.
$$

%%%%%%%%%%%%%%%%%%%%%%%%%%%%%%%%%%%%%%%%%
\section{Derivation of integral formula for trace}
The calculation of the trace using the bosonic expression of the
intertwining operators are exactly similar to the case of 
$sl_2$ \cite{JMMN}\cite{JM}.
Thus we present here only the necessary information
for the calculation of the trace.

We use the following formula
\begin{eqnarray}
&&
\hbox{te}_{V(\La_i)}\Big(
x^D
\exp\big(
\sum_{j=1}^{N-1}\sum_{n=1}^\infty
A^{(j)}_na_J(-n)
\big)
\exp\big(
\sum_{j=1}^{N-1}\sum_{n=1}^\infty
B^{(j)}_na_J(n)
\big)
\nonumber
\\
&&
\times
\exp(\sum_{j=1}^{N-1}c_j \alpha_j)
\prod_{j=1}^{N-1}g_j^{\partial_{\alpha_j}}
g_0^{-\partial_{\La_1}}
g_N^{-\partial_{\La_{N-1}}}
\Big)
\nonumber
\\
&&
=(x^N)_\infty^{-1}
g_0^{-(\La_1\vert\La_i)}
g_N^{-(\La_{N-1}\vert\La_i)}
g_i^{1-\delta_{i0}}
\theta_i(
g_0^{-1}g_1^2g_2^{-1},
\cdots,
g_{N-2}^{-1}g_{N-1}^2g_N^{-1}\vert x^N)
\nonumber
\\
&&
\times
\prod_{j=1}^{N-1}
\exp\Big(
\sum_{m=1}^\infty
\sum_{n=1}^\infty
{1\over n}
x^{Nmn}
A^{(j)}_n
(-[n]^2B^{(j-1)}_n+[n][2n]B^{(j)}_n-[n]^2B^{(j+1)}_n)
\Big),
\label{tracef}
\end{eqnarray}
where we set $B^{(0)}_n=B^{(N)}_n=0$.
The derivation of this formula is similar to
the $sl_2$ case.
We refer to \cite{JM} for details.

In the previous section we have given the expression of the
operators in terms of their normally ordered operators.
Therefore in this section we shall give a list of contributions
to the trace from the normally ordered operators.
Then using the formula (\ref{tracef}) we can calculate the
trace and the result is presented in section \ref{trf}.

For an operator ${\cal O}$ such that
\begin{eqnarray}
&&
{\cal O}=\exp\big(
\sum_{j=1}^{N-1}\sum_{n=1}^\infty
A^{(j)}_na_J(-n)
\big)
\exp\big(
\sum_{j=1}^{N-1}\sum_{n=1}^\infty
B^{(j)}_na_J(n)
\big)
\nonumber
\\
&&
\times
\exp(\sum_{j=1}^{N-1}c_j \alpha_j)
\prod_{j=1}^{N-1}g_j^{\partial_{\alpha_j}}
g_0^{-\partial_{\La_1}}
g_N^{-\partial_{\La_{N-1}}},
\nonumber
\end{eqnarray}
if we write
\begin{eqnarray}
&&
{\cal O}\approx J
\nonumber
\end{eqnarray}
then it means that 
\begin{eqnarray}
&&
J=
\prod_{j=1}^{N-1}
\exp\Big(
\sum_{m=1}^\infty
\sum_{n=1}^\infty
{1\over n}
x^{Nmn}
A^{(j)}_n
(-[n]^2B^{(j-1)}_n+[n][2n]B^{(j)}_n-[n]^2B^{(j+1)}_n)
\Big).
\nonumber
\end{eqnarray}

The following is the list which is necessary for the 
calculation of the trace.

$$
\tilde{\Phi}^{h(i)}_{N-1}(z)
\approx
{\{q^2x^N\} \over \{q^{2N}x^N\}},
\quad
\tilde{\Psi}^{\ast h(i)}_{N-1}(u)
\approx
{\{x^N\} \over \{q^{2N-2}x^N\}},
$$

$$
X^{-}_j(w)
\approx
(q^2x^N)_\infty 
(x^N)_\infty,
\quad
X^{+}_j(w)
\approx
(q^{-2}x^N)_\infty 
(x^N)_\infty 
$$

$$
:X^{\sigma_1}_{j_1}(w_1)X^{\sigma_2}_{j_2}(w_2):
\approx
1
\quad
\vert j_1-j_2\vert>1, \sigma_1,\sigma_2\in\{\pm\},
$$

$$
:X^{-}_{j}(w_1)X^{-}_{j+1}(w_2):
\approx
{
1
\over
(qx^Nw_1/w_2)_\infty
(qx^Nw_2/w_1)_\infty
},
$$

$$
:X^{-}_{j}(w_1)X^{-}_{j}(w_2):
\approx
(q^2x^Nw_1/w_2)_\infty
(q^2x^Nw_2/w_1)_\infty
(x^Nw_1/w_2)_\infty
(x^Nw_2/w_1)_\infty,
$$

$$
:X^{-}_{j}(w)X^{+}_{j+1}(v):
\approx
(x^Nw/v)_\infty
(x^Nv/w)_\infty,
$$

$$
:X^{-}_{j+1}(w)X^{+}_{j}(v):
\approx
(x^Nw/v)_\infty
(x^Nv/w)_\infty,
$$

$$
:X^{-}_{j}(w)X^{+}_{j}(v):
\approx
{
1
\over
(qx^Nw/v)_\infty
(qx^Nv/w)_\infty
(q^{-1}x^Nw/v)_\infty
(q^{-1}x^Nv/w)_\infty,
}
$$

$$
:X^{+}_{j}(v_1)X^{+}_{j+1}(v_2):
\approx
{
1
\over
(q^{-1}x^Nv_1/v_2)_\infty
(q^{-1}x^Nv_2/v_1)_\infty,
}
$$

$$
:X^{+}_{j}(v_1)X^{+}_{j}(v_2):
\approx
(x^N v_1/v_2)_\infty
(x^N v_2/v_1)_\infty
(q^{-2}x^N v_1/v_2)_\infty
(q^{-2}x^N v_2/v_1)_\infty,
$$

$$
:\tilde{\Phi}^{h(i)}_{N-1}(z)X^{\pm}_j(w):
\approx
1,
\quad
j\neq N-1,
$$

$$
:\tilde{\Psi}^{\ast h(i)}_{N-1}(u)X^{\pm}_j(w):
\approx
1,
\quad
j\neq N-1,
$$

$$
:\tilde{\Phi}^{h(i)}_{N-1}(z)X^{-}_{N-1}(w):
\approx
{
1
\over
(q^{N+2}x^Nz/w)_\infty
(q^{-N}x^N w/z)_\infty
},
$$

$$
:\tilde{\Phi}^{h(i)}_{N-1}(z)X^{+}_{N-1}(v):
\approx
(q^{N+1}x^Nz/v)_\infty
(q^{-N-1}x^N v/z)_\infty,
$$

$$
:\tilde{\Psi}^{\ast h(i)}_{N-1}(u)X^{-}_{N-1}(w):
\approx
(q^{N+1}x^N u/w)_\infty
(q^{-N-1}x^N w/u)_\infty,
$$

$$
:\tilde{\Psi}^{\ast h(i)}_{N-1}(u)X^{+}_{N-1}(v):
\approx
{
1
\over
(q^{N}x^N u/v)_\infty
(q^{-N-2}x^N v/u)_\infty
},
$$

$$
:\tilde{\Phi}^{h(i_1)}_{N-1}(z_1)
\tilde{\Phi}^{h(i_2)}_{N-1}(z_2):
\approx
{
\{q^2x^Nz_1/z_2\}
\{q^2x^Nz_2/z_1\}
\over
\{q^{2N}x^Nz_1/z_2\}
\{q^{2N}x^Nz_2/z_1\}
},
$$

$$
:\tilde{\Phi}^{h(i_1)}_{N-1}(z)
\tilde{\Psi}^{\ast h(i_2)}_{N-1}(u):
\approx
{
\{q^{2N-1}x^N z/u\}
\{q^{2N-1}x^N u/z\}
\over
\{qx^N z/u\}
\{qx^N u/z\}
},
$$

$$
:\tilde{\Psi}^{\ast h(i_1)}_{N-1}(u_1)
\tilde{\Psi}^{\ast h(i_2)}_{N-1}(u_2):
\approx
{
\{x^N u_1/u_2\}
\{x^N u_2/u_1\}
\over
\{q^{2N-2}x^N u_1/u_2\}
\{q^{2N-2}x^N u_2/u_1\}
}.
$$


\begin{thebibliography}{99}

\bibitem{BBS} Babelon, O., Bernard, D. and Smirnov, F.,
{Null-Vectors in Integrable Field Theory,
Commun. Math. Phys., {\bf 186} (1997) 601-648.}

\bibitem{DO} Date, E. and Okado, M.,
{Calculation of excitation spectra of the spin model
related with the vector representation of
the quantized affine algebra of type $A^{(1)}_n$,
Int. J. Mod. Phys. A, {\bf 9} (1994) 399-417.}

\bibitem{E} Etingof, P.,
{Difference equations with elliptic coefficients
and quantum affine algebras,
hep-th/9312057}
%Proc. Natl. Acad. Sci. USA, {\bf 85} (1988) 9373-9377.}

\bibitem{FJ} Frenkel, I. and Jing, N.,
{Vertex representations of quantum affine algebras,
Proc. Natl. Acad. Sci. USA, {\bf 85} (1988) 9373-9377.}

\bibitem{FK} Frenkel, I and Kac, V.,
{Basic representations of affine Lie algebras and
dual resonance models,
Invent. Math., {\bf 62} (1980) 23-66.}

\bibitem{FR} Frenkel, I. and Reshetikhin, N.,
{Quantum affine algebras and holonomic difference equations,
Commun. Math. Phys., {\bf 146} (1992) 1-60.}

\bibitem{HKMW} Hong, J. Kang, S-J., Miwa, T. and Weston, R.,
{Vertex models with alternating spins,
math.QA/9811175 {\bf }}

\bibitem{JM} Jimbo, M. and Miwa, T.,
{Algebraic analysis of solvable lattice models,
CBMS Regional Conference Series in Math. AMS, {\bf 85}, (1995).}

\bibitem{JMMN} Jimbo, M., Miki, K., Miwa, T. and Nakayashiki, A.,
{Correlation function of the XXZ model for $\Delta<-1$,
Phys. Lett. A, {\bf 168} (1992) 256-263.}

\bibitem{Ko} Koyama, Y.,
{Staggered polarization of vertex models with
$\uqn$ symmetry,
Commun. Math. Phys., {\bf 164} (1994) 277-291.}

\bibitem{Na} Nakayashiki, A.,
{Fusion of $q$-vertex operators and its application to solvable
vertex models,
Commun. Math. Phys., {\bf 177} (1996) 27-62.}

\bibitem{NPT} Nakayashiki, A., Pakuliak, S. and Tarasov, V.,
{On solutions of the KZ and qKZ equations at level zero,
to appear in Ann. Inst. Henri Poincare, q-alg/9712002.}

\bibitem{S} Smirnov, F.,
{Counting the local fields in SG theory,
Nucl. Phys. B {\bf 453} (1995) 807-824.}

\bibitem{T} Tarasov, V.,
{Completeness of the hypergeometric solutions of the qKZ equation
at level zero,
Max-Planck-Inst. Preprint Series 1998 (87).}

\bibitem{TV} Tarasov, V. and Varchenko, A.,
{Geometry of q-hypergeometric functions, quantum affine algebras,
and elliptic quantum groups,
Asterisque {\bf 246} (1997) 1-135.}

\end{thebibliography}
\end{document}